\documentclass{agujournal2019}

\usepackage{algorithm}
\usepackage{algpseudocode}
\usepackage{amsfonts}
\usepackage{makecell}
\usepackage{mathtools}
\usepackage{multicol}
\usepackage{multirow}
\usepackage{relsize}

\hypersetup{
	colorlinks=true,
	linkcolor=black,
	citecolor=black,
	filecolor=black,
	urlcolor=black,
}

\newcommand{\doi}[1]{\url{https://doi.org/#1}}

\begin{document}
	\title{Optimization of Structural Flood Mitigation Strategies}

	\authors{Byron Tasseff\affil{1, 2}, Russell Bent\affil{3}, and Pascal Van Hentenryck\affil{4}}
	\affiliation{1}{Department of Industrial and Operations Engineering, University of Michigan, Ann Arbor, Michigan, USA.}
	\affiliation{2}{Information Systems and Modeling Group, Los Alamos National Laboratory, Los Alamos, New Mexico, USA.}
	\affiliation{3}{Applied Mathematics and Plasma Physics Group, Los Alamos National Laboratory, Los Alamos, New Mexico, USA.}
	\affiliation{4}{H. Milton Stewart School of Industrial and Systems Engineering, Georgia Institute of Technology, Atlanta, Georgia, USA.}
	\correspondingauthor{Byron Tasseff}{btasseff@lanl.gov}

	\begin{keypoints}
		\item The structural optimal flood mitigation problem is introduced
		\item A problem discretization amenable to derivative-free optimization is developed
		\item Benefits of constraining the problem with additional physics-based restrictions are shown
	\end{keypoints}

	\begin{abstract}
		The dynamics of flooding are primarily influenced by the shape, height, and roughness (friction) of the underlying topography.
		For this reason, mechanisms to mitigate floods frequently employ structural measures that either modify topographic elevation, e.g., through the placement of levees and sandbags, or increase roughness, e.g., through revegetation projects.
		However, the configuration of these measures is typically decided in an ad hoc manner, limiting their overall effectiveness.
		The advent of high-performance surface water modeling software and improvements in black-box optimization suggest that a more principled design methodology may be possible.
		This paper proposes a new computational approach to the problem of designing structural mitigation strategies under physical and budgetary constraints.
		It presents the development of a problem discretization amenable to simulation-based, derivative-free optimization.
		However, meta-heuristics alone are found to be insufficient for obtaining quality solutions in a reasonable amount of time.
		As a result, this paper proposes novel numerical and physics-based procedures to improve convergence to a high-quality mitigation.
		The efficiency of the approach is demonstrated on hypothetical dam break scenarios of varying complexity under various mitigation budget constraints.
		In particular, experimental results show that, on average, the final proposed algorithm results in a $65\%$ improvement in solution quality compared to a direct implementation.
	\end{abstract}

	\section{Introduction}
	\label{sec:introduction}
	Modern flood risk management (FRM) is a continuous process of identifying issues, defining objectives, assessing risks, appraising options, implementation, monitoring, and review.
	Within this framework, risk assessment is regarded as a cyclic process that includes the design and evaluation of alternative management strategies.
	Such strategies commonly include both ``hard'' and ``soft'' structural mitigation measures, e.g., the construction of dams (hard) and wetland storage (soft) \cite{sayers2013flood}.
	Measures can also be temporary (e.g., sandbags) or permanent (e.g., levees).
	However, for complex scenarios, the number of feasible strategies is extremely large and computationally difficult to explore.
	As such, the manual design and assessment of these strategies, whether conducted in a real-world or simulation-based setting, can be time-consuming and expensive.
	This limitation may result in vastly suboptimal FRM strategies.
	To aid in the FRM process, an optimization-based decision support approach for proposing structural mitigation designs can serve as a useful tool within the overall risk assessment phase.

	This paper develops such an optimization-based decision support approach for proposing flood protection strategies, whereby effective mitigation designs are realized through the exploration of various configurations in a computational setting.
	Specifically, the paper defines the Optimal Flood Mitigation Problem (OFMP), whose goal is to make topographic modifications that protect critical regions under a given flood scenario.
	This is a difficult optimization problem, as these modifications can have highly nonlinear effects on the flooding behavior.
	Moreover, physical models used to examine these effects are computationally expensive.
	Finally, as the OFMP aims at deciding several simultaneous modifications, an efficient exploration of the full search space is computationally intractable for realistic scenarios.

	The literature associated with the OFMP is limited.
	The closest related studies are by \citeA{jtbp14flood} and \citeA{tasseff2016optimal}.
	In the former, an interdiction model for flood mitigation is proposed, and model surrogates constructed from simulation data are used as proxies for estimating flood sensitivity to hard structural mitigation measures.
	In the latter, an OFMP similar to that discussed herein is introduced, and mixed-integer linear programs constrained by approximate flooding dynamics are solved to obtain hard structural mitigation designs.
	However, the approach is shown to suffer from substantial scalability issues in space and time.
	Neither study uses an approach which relies upon repeated deterministic modeling of the partial differential equations (PDEs) underlying the flooding dynamics.

	A number of studies discuss simulation-optimization approaches for reservoir operation, where the PDEs associated with the flood dynamics are treated as a black box.
	An extensive literature review of these studies can be found in \citeA{Che2015}.
	The work described by \citeA{Colombo2009} considers the full PDEs, but their focus is on optimizing normal operations of an open-channel system.
	Finally, the problem of optimizing dike heights with uncertainty in flooding estimates is considered by \citeA{Brekelmans2012}.
	However, in this study, the PDEs for flood propagation are not considered, and probability models for maximum flood depths are used in place of deterministic physical models.

	This paper presents a new approach to the problem of designing structural FRM strategies over PDE constraints.
	It develops a problem discretization amenable to simulation-based derivative-free optimization.
	Moreover, the paper shows that meta-heuristics alone are insufficient for obtaining quality solutions in reasonable time.
	As a result, it presents several innovative computational and physics-based techniques to increase convergence to high-quality solutions.
	The efficiency of the proposed approach is compared using hypothetical dam break scenarios of varying complexity under multiple mitigation budgets.
	Experimental results show that the proposed algorithm results in a $65\%$ improvement in solution quality compared to a direct implementation.

	The rest of this paper is organized as follows:
	Section \ref{sec:model} discusses the background of flood modeling and formalization of the OFMP;
	Section \ref{sec:algorithm} describes solution methods for a specific OFMP;
	Section \ref{sec:results} compares these methods using fictional dam break scenarios, with both simplistic (Section \ref{subsect:results-simple}) and realistic (Section \ref{subsect:results-icold}) topographies, and multiple mitigation budgets;
	and Section \ref{sec:conclusion} concludes the paper.

	\section{Model}
	\label{sec:model}

	In this paper, it is assumed that flood scenarios are modeled using the two-dimensional (2D) shallow water equations.
	These PDEs are derived from the Navier-Stokes equations under the assumption that horizontal length scales are much larger than the vertical scale.
	This is reasonable for large-scale floods, where water depths are much smaller than typical flood wavelengths.
	Two-dimensional models, in particular, alleviate the fundamental disadvantages of their 1D counterparts by allowing for higher-order representations of the topographic surface.
	Moreover, 2D models readily make use of widely available topographic elevation data.
	Finally, with recent advances in high-performance computing, solutions to these PDEs have become numerically tractable for large-scale problems, making them of particular computational interest.
	With volumetric, bed slope, and bed shear stress source terms, these equations are expressed as
	\begin{linenomath*}
	\begin{subequations}
	\begin{align}
		\frac{\partial h}{\partial t} + \frac{\partial (h u)}{\partial x} + \frac{\partial (h v)}{\partial y} &= R(x, y, t), \label{eqn:2d-swe-1} \\
		\frac{\partial (h u)}{\partial t}+ \frac{\partial}{\partial x}\left( h u^2 + \frac{1}{2}g h^2 \right) + \frac{\partial (h u v)}{\partial y} &= -g h \frac{\partial B}{\partial x} - g \frac{n^{2}}{h^{1/3}} |u| u, \label{eqn:2d-swe-2} \\
		\frac{\partial (h v)}{\partial t} + \frac{\partial (h uv)}{\partial x} + \frac{\partial}{\partial y}\left(h v^2 + \frac{1}{2}g h ^2\right) &= -g h \frac{\partial B}{\partial y} - g \frac{n^{2}}{h^{1/3}} |v| v, \label{eqn:2d-swe-3}
	\end{align}
	\end{subequations}
	\end{linenomath*}
	where $h$ is the water depth, $u$ and $v$ are horizontal velocities, $B$ is the bottom topography (or bathymetry), $g$ is the acceleration due to gravity, $n$ is the Manning's roughness coefficient, and $R$ is a volumetric source term \cite{FLD:FLD4023}.
	Equation \eqref{eqn:2d-swe-1} represents mass continuity, while Equations \eqref{eqn:2d-swe-2} and \eqref{eqn:2d-swe-3} represent conservation of momentum over the two horizontal dimensions.

	These equations can be rewritten in vector form by introducing the definitions
	\begin{linenomath*}
	\begin{equation}
	\begin{gathered}
	\mathbf{U} := \left(h, hu, hv\right), ~ \mathbf{F}(\mathbf{U}) := \left(hu, hu^{2} + \frac{1}{2} g h^{2}, huv\right), \\
	\mathbf{G}(\mathbf{U}) := \left(hv, huv, hv^{2} + \frac{1}{2} g h^{2}\right), ~ \mathbf{S}_{R}(R) := \left(R(x, y, t), 0, 0\right), \\
	\mathbf{S}_{B}(\mathbf{U}, B) := \left(0, -g h \frac{\partial B}{\partial x}, -g h \frac{\partial B}{\partial y}\right), ~ \mathbf{S}_{n} := \left(0, -g \frac{n^{2}}{h^{1/3}} |u| u, -g \frac{n^{2}}{h^{1/3}} |v| v\right),
	\end{gathered}
	\end{equation}
	\end{linenomath*}
	where $\mathbf{U}$ is the vector of conserved variables; $\mathbf{F}$ and $\mathbf{G}$ are fluxes in the $x$- and $y$-directions, respectively; and $\mathbf{S}_{R}$, $\mathbf{S}_{B}$, and $\mathbf{S}_{n}$ are the volumetric, bed slope, and bed shear stress source terms, respectively.
	This allows Equations \eqref{eqn:2d-swe-1}, \eqref{eqn:2d-swe-2}, and \eqref{eqn:2d-swe-3} to be rewritten more concisely as
	\begin{linenomath*}
	\begin{equation}
	\mathbf{U}_{t} + \mathbf{F}_{x} + \mathbf{G}_{y} = \mathbf{S}_{R} + \mathbf{S}_{B} + \mathbf{S}_{n},
	\end{equation}
	\end{linenomath*}
	where $t$, $x$, and $y$ indicate partial differentiation with respect to those variables.

	The OFMP considers a flood scenario (e.g., a dam failure) and a set of 2D regions (``assets'') to protect.
	To minimize flooding at asset locations, the model must produce optimal topographic elevation and roughness fields using a set of $m$ mitigation measures.
	For each measure $i \in \{1, 2, \dots, m\}$, the functions $\delta_{B}(\omega_{i})$ and $\delta_{n}(\omega_{i})$ define 2D fields of height and roughness for a given tuple of field parameters $\omega_{i}$.
	Measures can first additively modify the elevation field $B$ to return a new field $\tilde{B}$, defined as
	\begin{linenomath*}
	\begin{equation}
		\tilde{B}(B, \left(\omega_{1}, \omega_{2}, \dots, \omega_{m}\right)) := B + \sum_{i = 1}^{m} \delta_{B}(\omega_{i}) \label{eqn:elevation-configuration}.
	\end{equation}
	\end{linenomath*}
	Similarly, measures can modify the roughness field $n$ to return a new field defined as
	\begin{linenomath*}
	\begin{equation}
		\tilde{n}(n, \left(\omega_{1}, \omega_{2}, \dots, \omega_{m}\right)) := n + \{\max_{i}\{\delta_{n}(\omega(i))(x, y)\} : (x, y) \in \mathbb{R}^{2}\} \label{eqn:roughness-configuration},
	\end{equation}
	\end{linenomath*}
	i.e., a field of maximum roughness.
	For notational ease, hereafter, $\tilde{B}$ refers to Equation \ref{eqn:elevation-configuration}, $\tilde{n}$ refers to Equation \ref{eqn:roughness-configuration}, and the tuple $(\omega_{1}, \omega_{2}, \dots, \omega_{m})$ is referred to as the ``parametric configuration,'' or simply the ``configuration.''
	With these definitions and shorthand notations, the \emph{modified} bed slope source term is defined as
	\begin{linenomath*}
	\begin{equation}
	\tilde{\mathbf{S}}_{B}\left(\mathbf{U}, B, \left(\omega_{1}, \omega_{2}, \dots, \omega_{m}\right)\right) := \left(0, -g h \frac{\partial \tilde{B}}{\partial x}, -g h \frac{\partial \tilde{B}}{\partial y}\right).
	\end{equation}
	\end{linenomath*}
	We note that the change in elevation may be a result of permanent structures such as levees or temporary measures such as sandbags.
	Similarly, the modified bed shear stress source term is defined as
	\begin{linenomath*}
	\begin{equation}
	\tilde{\mathbf{S}}_{n}\left(\mathbf{U}, n, \left(\omega_{1}, \omega_{2}, \dots, \omega_{m}\right)\right) := \left(0, -g \frac{\tilde{n}^{2}}{h^{1/3}} |u| u, -g \frac{\tilde{n}^{2}}{h^{1/3}} |v| v\right).
	\end{equation}
	\end{linenomath*}
	Hereafter, $\tilde{\mathbf{S}}_{B} := \tilde{\mathbf{S}}_{B}\left(\mathbf{U}, B, (\omega_{1}, \omega_{2}, \dots, \omega_{m})\right)$ and $\tilde{\mathbf{S}}_{n} := \tilde{\mathbf{S}}_{n}\left(\mathbf{U}, n, (\omega_{1}, \omega_{2}, \dots, \omega_{m})\right)$ are used to concisely denote these two source terms that vary with the configuration.

	The OFMP is then written in a form that embeds the 2D shallow water equations as constraints and optimizes the tuple $(\omega_{1}, \omega_{2}, \dots, \omega_{m})$ (i.e., the configuration) via
	\begin{linenomath*}
	\begin{subequations}
	\begin{align}
		 & \underset{\omega_{1}, \omega_{2}, \dots, \omega_{m}}{\text{minimize}}
		 & & z\left(\omega_{1}, \omega_{2}, \dots, \omega_{m}\right) = \mathlarger{\sum}_{a \in \mathcal{A}} \iint_{a} \max_{t} {h(x, y, t)} \,dx \,dy \label{eqn:general-ofmp-objective} \\
		 & \text{subject to}
		 & & \mathbf{U}_{t} + \mathbf{F}_{x} + \mathbf{G}_{y} = \mathbf{S}_{R} + \tilde{\mathbf{S}}_{B} + \tilde{\mathbf{S}}_{n} \label{eqn:general-ofmp-2d-swe} \\
		 & & & \delta_{B}(\omega_{i})(x, y) = 0, ~ \forall i \in \{1, 2, \dots, m\}, ~ \text{for} ~ (x, y) \in \bigcup \mathcal{A} \label{eqn:general-ofmp-prohibit-structures} \\
		 & & & \delta_{n}(\omega_{i})(x, y) = 0, ~ \forall i \in \{1, 2, \dots, m\}, ~ \text{for} ~ (x, y) \in \bigcup \mathcal{A} \label{eqn:general-ofmp-prohibit-roughness} \\
		 & & & \left(\omega_{1}, \omega_{2}, \dots, \omega_{m}\right) \in \mathcal{F} \label{eqn:general-ofmp-valid-structures}.
	\end{align}
	\end{subequations}
	\end{linenomath*}
	Here, $\mathcal{A}$ denotes the set of asset regions to be protected and $z$ denotes the objective function.
	This function is defined in Equation \eqref{eqn:general-ofmp-objective} and captures the maximum water volume over all asset locations and times.
	Constraint \eqref{eqn:general-ofmp-2d-swe} denotes the solution to the shallow water equations in the presence of the $m$ mitigation measures.
	Constraints \eqref{eqn:general-ofmp-prohibit-structures} prohibit measures from being constructed ``underneath'' an asset.
	Similarly, Constraints \eqref{eqn:general-ofmp-prohibit-roughness} prohibit the roughness at an asset location from being modified.
	Finally, Constraint \eqref{eqn:general-ofmp-valid-structures} ensures $(\omega_{1}, \omega_{2}, \dots, \omega_{m})$ resides within the set of all feasible parametric configurations $\mathcal{F}$, i.e., $\mathcal{F}$ distinguishes valid and invalid mitigation designs.

	For simplicity of presentation, this paper considers only two types of structural measures, although the approach can easily be generalized to include other soft and hard measures, both temporary and permanent.
	The first type is an immovable wall of fixed length ($\ell$), width ($w$), and height ($\bar{b}_{i}$).
	Each wall is defined using three continuously-defined, bounded parameters: latitudinal position of the wall centroid ($\lambda_{i}$), longitudinal position of the wall centroid ($\phi_{i}$), and angle of the wall formed with respect to the longitudinal axis ($\theta_{i}$).
	In this paper, the centroid position is bounded by the scenario domain's spatial extent, and $\theta_{i} \in [0, \pi]$.
	The second structural type is a revegetation project defined by a 2D circular region with center $(\lambda_{i}, \phi_{i})$ and fixed radius $r$ that increases the area's Manning's roughness coefficient based on a fixed field $\bar{n}_{i}$.
	Under these assumptions, the OFMP aims at deciding $\omega_{i} = \left(\lambda_{i}, \phi_{i}, \theta_{i}, \bar{b}_{i}, \bar{n}_{i}\right)$ for each measure $i \in \{1, 2, \dots, m\} = \mathcal{M}$, where $\bar{b}_{i}$ and $\bar{n}_{i}$ are decided a priori for each measure.
	More specifically, this produces an OFMP of the specialized form
	\begin{linenomath*}
	\begin{subequations}
	\begin{align}
		 & \underset{\omega_{1}, \omega_{2}, \dots, \omega_{m}}{\text{minimize}}
		 & & z\left(\omega_{1}, \omega_{2}, \dots, \omega_{m}\right) = \mathlarger{\sum}_{a \in \mathcal{A}} \iint_{a} \max_{t} {h(x, y, t)} \,dx \,dy \label{eqn:specific-ofmp-objective} \\
		 & \text{subject to}
		 & & \mathbf{U}_{t} + \mathbf{F}_{x} + \mathbf{G}_{y} = \mathbf{S}_{R} + \tilde{\mathbf{S}}_{B} + \tilde{\mathbf{S}}_{n} \label{eqn:specific-ofmp-2d-swe} \\
		 & & & \delta_{B}(\omega_{i})(x, y) = 0, ~ \text{for} ~ (x, y) \in \bigcup \mathcal{A}, ~ \forall i \in \mathcal{M} \label{eqn:specific-ofmp-no-overlap} \\
		 & & & \delta_{n}(\omega_{i})(x, y) = 0, ~ \text{for} ~ (x, y) \in \bigcup \mathcal{A}, ~ \forall i \in \mathcal{M} \label{eqn:specific-ofmp-no-manning} \\
		 & & & \delta_{B}(\omega_{i})(x, y) =
			  \begin{cases}
			  \begin{aligned}
				  \bar{b}_{i} & ~ \textrm{for}
						 \begin{cases}
							  |(x - \phi_{i}) \cos \theta_{i} - (y - \lambda_{i}) \sin \theta_{i}| \leq \frac{\ell}{2} \\
							  |(x - \phi_{i}) \sin \theta_{i} + (y - \lambda_{i}) \cos \theta_{i}| \leq \frac{w}{2}
						 \end{cases} \\
			  0 & ~ \textrm{otherwise}
			  \end{aligned}
		  \end{cases} \hspace{-2.05em} \forall i \in \mathcal{M} \label{eqn:specific-ofmp-wall-definition} \\
		 & & & \delta_{n}(\omega_{i})(x, y) =
			  \begin{cases}
				  \bar{n}_{i}(x, y) & \textrm{for} ~ (x - \phi_{i})^{2} + (y - \lambda_{i})^{2} \leq r^{2} \\
				  0 & \textrm{otherwise}
			  \end{cases} \forall i \in \mathcal{M} \label{eqn:specific-ofmp-wall-definition-veg} \\
		 & & & \lambda_{lb} \leq \lambda_{i} \leq \lambda_{ub}, ~ \phi_{lb} \leq \phi_{i} \leq \phi_{ub}, ~ 0 \leq \theta_{i} \leq \pi, ~ \forall i \in \mathcal{M} \label{eqn:specific-ofmp-valid-structures}.
	\end{align}
	\end{subequations}
	\end{linenomath*}
	Using this formulation, $i$ is a wall when $\bar{b}_{i} > 0$ and $\bar{n}_{i} = 0$, and $i$ is a revegetation project when $\bar{b}_{i} = 0$ and $\bar{n}_{i} > 0$.
	Constraints \eqref{eqn:specific-ofmp-no-overlap} and \eqref{eqn:specific-ofmp-no-manning} emphasize that modifications cannot be made within asset regions; Constraints \eqref{eqn:specific-ofmp-wall-definition} impose the wall height $\bar{b}_{i}$ within each rotated rectangle defined using the parameters $\lambda_{i}$, $\phi_{i}$, and $\theta_{i}$ and a standard 2D rotation matrix; and Constraints \eqref{eqn:specific-ofmp-wall-definition-veg} impose additions to roughness within each revegetation circle defined by the center $(\lambda_{i}, \phi_{i})$.
	Finally, Constraints \eqref{eqn:specific-ofmp-valid-structures} replace Constraint \eqref{eqn:general-ofmp-valid-structures} of the more general OFMP.
	Here, $\lambda_{lb}$ and $\lambda_{ub}$ ($\phi_{lb}$ and $\phi_{ub}$) are the lower and upper latitudinal (longitudinal) boundaries of the scenario domain.

	Constraints \eqref{eqn:specific-ofmp-valid-structures} imply a large feasible region, as the spatial extent is typically much larger than the flood's extent.
	To reduce the solution space, the notion of a \emph{restricted region} $\mathcal{P}$ is thus introduced, where centroids must reside in $\mathcal{P}$.
	That is,
	\begin{linenomath*}
	\begin{equation}
	(\lambda_{i}, \phi_{i}) \in \mathcal{P}, ~ \forall i \in \{1, 2, \dots, m\} \tag{9h} \label{eqn:specific-ofmp-centroid-restriction}
	\end{equation}
	\end{linenomath*}
	is appended to the problem above, completing the primary model used in this paper.

	\section{Algorithm}
	\label{sec:algorithm}
	The OFMP at the end of Section \ref{sec:model} remains difficult to solve directly.
	However, with recent improvements in both numerical discretizations of the shallow water equations (e.g., \citeA{FLD:FLD4023}) and high-performance implementations thereof (e.g., \citeA{brodtkorb2012efficient}, \citeA{tasseff2016nuflood}), numerically efficient solutions of the PDEs described in Constraint \eqref{eqn:specific-ofmp-2d-swe} are possible.
	With this intuition, in Algorithm \ref{alg:pseudoalg}, a time-limited search-based method is introduced to find a near-optimal solution $\left(\omega_{1}^{*}, \omega_{2}^{*}, \dots, \omega_{m}^{*}\right)$ to the problem defined by Equations \eqref{eqn:specific-ofmp-objective} through \eqref{eqn:specific-ofmp-centroid-restriction}.

	\begin{linenomath*}
	\begin{algorithm}[ht]
		\caption{\textproc{SolveOFMP}: Solves the OFMP defined by Equations \eqref{eqn:specific-ofmp-objective} through \eqref{eqn:specific-ofmp-centroid-restriction}.}
		 \label{alg:pseudoalg}
		 \begin{algorithmic}[1]
		 \Function{SolveOFMP}{$B, n, \mathcal{A}, m, T_{\max}, \alpha$}
			  \State $\tilde{\mathcal{P}} \gets \textproc{InitializeRestriction}(B, n, \mathcal{A}, \alpha)$ \label{alg:pseudoalg:line:init-2}
			  \State $\left(\omega_{1}^{*}, \omega_{2}^{*}, \dots, \omega_{m}^{*}\right) \gets \textproc{InitializeSolution}(m, \tilde{\mathcal{P}}), ~ \Omega \gets \emptyset$ \label{alg:pseudoalg:line:init-1}
			  \While{$\textproc{Clock} < T_{\max}$}
					\State $\left(\omega_{1}, \omega_{2}, \dots, \omega_{m}\right) \gets \textproc{GenerateSolution}(m, \tilde{\mathcal{P}}, \Omega)$ \label{alg:pseudoalg:line:generate-solution}
					\State $\textrm{Solve} ~ \mathbf{U}_{t} + \mathbf{F}_{x} + \mathbf{G}_{y} = \mathbf{S}_{R} + \tilde{\mathbf{S}}_{B} + \tilde{\mathbf{S}}_{n}$ \label{alg:pseudoalg:line:solve-flood-pdes}
					\State $\Omega \gets \Omega \cup \left\{\left(\omega_{1}, \omega_{2}, \dots, \omega_{m}\right)\right\}$ \label{alg:pseudoalg:line:update-solution-history}
					\If{$z(\omega_{1}, \omega_{2}, \dots, \omega_{m}) < z(\omega^{*}_{1}, \omega^{*}_{2}, \dots, \omega^{*}_{m})$} \label{alg:pseudoalg:line:start-if-update-best-solution}
						 \State $\left(\omega_{1}^{*}, \omega_{2}^{*}, \dots, \omega_{m}^{*}\right) \gets \left(\omega_{1}, \omega_{2}, \dots, \omega_{m}\right)$ \label{alg:pseudoalg:line:update-best-solution}
						 \State $\tilde{\mathcal{P}} \gets \textproc{UpdateRestriction}(\mathbf{U}, \mathcal{A}, \tilde{\mathcal{P}}, \alpha)$ \label{alg:pseudoalg:line:update-restricted-region}
					\EndIf \label{alg:pseudoalg:line:end-if-update-best-solution}
			  \EndWhile
			  \State \Return $(B + \sum_{i = 1}^{m} \delta_{B}(\omega^{*}_{i})$, $n + \{\max_{i}\{\delta_{n}(\omega^{*}_{i})(x, y)\} : (x, y) \in \mathbb{R}^{2}\})$ \label{alg:pseudoalg:line:return-value}
			  \EndFunction
		 \end{algorithmic}
	\end{algorithm}
	\end{linenomath*}

	Here, $B$ and $n$ denote the initial topographic elevation and Manning's roughness coefficient fields; $\mathcal{A}$ denotes the set of assets; $m$ denotes the number of mitigation measures being configured; $T_{\max}$ denotes the maximum clock time; and $\alpha$ is a parameter used for computing restrictions.
	The function $\textproc{Clock}$ returns the current clock time.
	Since a useful definition of $\mathcal{P}$ is difficult to compute a priori, $\tilde{\mathcal{P}}$ serves as an iterative approximation of some desired $\mathcal{P}$.
	In Line \ref{alg:pseudoalg:line:init-2}, $\tilde{\mathcal{P}}$ is initialized; it is later modified in Line \ref{alg:pseudoalg:line:update-restricted-region} using $\textproc{UpdateRestriction}$.
	Both functions are described in Section \ref{sec:algorithm-initialize-update-restricted-region}.
	In Line \ref{alg:pseudoalg:line:init-1}, the best solution and the historical solution set $\Omega$ are initialized.
	In Line \ref{alg:pseudoalg:line:generate-solution}, a configuration is generated via some history-dependent function $\textproc{GenerateSolution}$, described in Section \ref{subsec:experimental_setting}.
	In Line \ref{alg:pseudoalg:line:solve-flood-pdes}, the shallow water equations are solved.
	In Line \ref{alg:pseudoalg:line:update-solution-history}, the historical solution set is updated.
	In Lines \ref{alg:pseudoalg:line:start-if-update-best-solution} through \ref{alg:pseudoalg:line:end-if-update-best-solution}, the best solution and $\tilde{\mathcal{P}}$ are updated.
	Finally, in Line \ref{alg:pseudoalg:line:return-value}, the best elevation and roughness fields are returned.

	\subsection{Computation of the restricted region}
	\label{sec:algorithm-initialize-update-restricted-region}
	\subsubsection{The direct methodology}
	The most obvious globally acceptable method for selecting $\mathcal{P}$ is to assume
	\begin{linenomath*}
	\begin{equation}
	\mathcal{P} = \mathbb{R}^{2},
	\end{equation}
	\end{linenomath*}
	where, of course,
	\begin{linenomath*}
	\begin{equation}
	\left\{(x, y) \in \mathbb{R}^{2} : \lambda_{lb} \leq x \leq \lambda_{ub}, ~ \phi_{lb} \leq y \leq \phi_{ub}\right\} \subset \mathcal{P},
	\end{equation}
	\end{linenomath*}
	indicating the bounds within Constraints \ref{eqn:specific-ofmp-valid-structures} involving $\lambda_{i}$ and $\phi_{i}$ dominate those imposed by $\mathcal{P}$.
	This method for selecting $\mathcal{P}$ is hereafter referred to as the direct method.
	In practice, this method is used to define the direct implementations of the functions \textproc{InitializeRestriction} and \textproc{UpdateRestriction}, both of which return the set $\mathbb{R}^{2}$.

	\subsubsection{The pathline methodology}
	A \emph{pathline} is the trajectory an individual fluid element follows over time, beginning at position $\left(x_{0}, y_{0}\right)$ and time $t_{0}$.
	In 2D, a pathline is defined by the two equations
	\begin{linenomath*}
	\begin{subequations}
	\begin{align}
	x(t) &= x_{0} + \int_{t_{0}}^{t} u(x(t^{\prime}), y(t^{\prime}), t^{\prime}) d t^{\prime}, \\
	y(t) &= y_{0} + \int_{t_{0}}^{t} v(x(t^{\prime}), y(t^{\prime}), t^{\prime}) d t^{\prime},
	\end{align}
	\end{subequations}
	\end{linenomath*}
	where $u$ and $v$ are velocities in the $x$- and $y$-directions.
	To compute the pathline from a flood wave to an initially dry point $(x_{0}, y_{0})$, the definition of $t_{\text{wet}}(x_{0}, y_{0})$ is introduced as the time at which the depth at $(x_{0}, y_{0})$ exceeds some threshold.
	More concisely,
	\begin{linenomath*}
	\begin{equation}
	t_{\text{wet}}(x_{0}, y_{0}) := \min\left\{t \in [t_{0}, t_{f}] : h(x_{0}, y_{0}, t) \geq \epsilon_{h}\right\} \label{eqn:continuous-wet-time},
	\end{equation}
	\end{linenomath*}
	where $\epsilon_{h}$ is an arbitrarily small depth, taken in this study to be one millimeter.
	Using this definition, the pathline equations may be integrated \emph{in reverse}, giving
	\begin{linenomath*}
	\begin{subequations}
	\begin{align}
	x_{\text{wet}}(x_{0}, y_{0}, t) &= x_{0} + \int_{t_{\text{wet}}}^{t} u(x_{\text{wet}}(t^{\prime}), y_{\text{wet}}(t^{\prime}), t^{\prime}) d t^{\prime}, \label{eqn:pathline-x} \\
	y_{\text{wet}}(x_{0}, y_{0}, t) &= y_{0} + \int_{t_{\text{wet}}}^{t} v(x_{\text{wet}}(t^{\prime}), y_{\text{wet}}(t^{\prime}), t^{\prime}) d t^{\prime}, \label{eqn:pathline-y}
	\end{align}
	\end{subequations}
	\end{linenomath*}
	where it is assumed that $t \leq t_{\text{wet}}$.
	The above equations approximate a path to flooding.

	In this paper, a \emph{pathtube} is defined as a set of pathlines satisfying Equations \eqref{eqn:pathline-x} and \eqref{eqn:pathline-y}.
	For a region $\mathcal{R}$, the pathtube $\mathcal{S}$ encompassing $\mathcal{R}$ with a start time of $t_{0}$ is
	\begin{linenomath*}
	\begin{equation}
	\mathcal{S}(\mathbf{U}, \mathcal{R}) = \left\{(x_{\text{wet}}(x_{0}, y_{0}, t), y_{\text{wet}}(x_{0}, y_{0}, t)) \in \mathbb{R}^{2} : (x_{0}, y_{0}) \in \mathcal{R}, ~ t \in [t_{0}, t_{\text{wet}}(x_{0}, y_{0})]\right\}.
	\end{equation}
	\end{linenomath*}
	This region encompasses approximate paths of least resistance from a flood to $\mathcal{R}$.
	It is clear that good locations for structural mitigation measures are likely to reside in $\mathcal{S}$.

	A robust selection of $\mathcal{P}$ would account for the change in $\mathbf{U}$ with respect to a large set of feasible configurations.
	In an ideal setting, a good selection for $\mathcal{P}$ would thus be
	\begin{linenomath*}
	\begin{equation}
	\mathcal{P} = \bigcup_{\boldsymbol{\omega} \in \mathcal{F}} \bigcup_{a \in \mathcal{A}} \left\{(x, y) \in \mathcal{S}(\mathbf{U}, a) : \mathbf{U}_{t} + \mathbf{F}_{x} + \mathbf{G}_{y} = \mathbf{S}_{R} + \tilde{\mathbf{S}}_{B} + \tilde{\mathbf{S}}_{n}\right\} \label{eqn:streamline-restricted-region}.
	\end{equation}
	\end{linenomath*}
	In practice, defining $\mathcal{P}$ as per Equation \eqref{eqn:streamline-restricted-region} is nontrivial.
	First, each $a \in \mathcal{A}$ may be a set of infinitely many points.
	There are also infinitely many moments $t$ in a solution $\mathbf{U}$ to the shallow water equations.
	Most importantly, the union over all feasible configurations $(\omega_{1}, \omega_{2}, \dots, \omega_{m}) = \boldsymbol{\omega} \in \mathcal{F}$ assumes knowledge of $\mathbf{U}$ for any such feasible configuration $(\omega_{1}, \omega_{2}, \dots, \omega_{m})$.
	For these reasons, an iteratively-constructed definition of the pathtube-like region $\tilde{\mathcal{P}}$ is instead proposed, which approximately captures the features of some unknown larger $\mathcal{P}$ relevant to the OFMP (e.g., Equation \eqref{eqn:streamline-restricted-region}).

	From a numerical perspective, each $a \in \mathcal{A}$ is actually a polygon whose exterior connects a set of points $P_{a}$.
	Solutions to the OFMP are likely to intersect the pathlines from a flood to each of these points.
	Also, in practice, numerical solutions to the shallow water equations are discrete in space and time.
	Assuming that solutions are obtained for a set of timestamps $\mathcal{T}$ on a rectangular grid $G$, $t_{\text{wet}}$ is first redefined as
	\begin{linenomath*}
	\begin{equation}
	t_{\text{wet}}(x_{0}, y_{0}) := \min\{t \in \mathcal{T} : h_{i_{0}, j_{0}, t} \geq \epsilon_{h}\} \label{eqn:discrete-wet-time},
	\end{equation}
	\end{linenomath*}
	where $(i_{0}, j_{0})$ is the unique index of the cell in grid $G$ that contains the point $(x_{0}, y_{0})$.

	For each point along an asset exterior, a numerical representation of the pathline leading \emph{to} that point is desired.
	To accomplish this, it is assumed that a pathline can be approximated as a set $\mathcal{L}$ of discrete points.
	These points can be generated by solving Equations \eqref{eqn:pathline-x} and \eqref{eqn:pathline-y} using any suitable ordinary differential equation integration technique.
	In this study, suggestions from \citeA{telea2014data} (initially described for \emph{streamlines}, which trace a static field) are used to compute pathlines according to the function $\textproc{ComputePathline}(\mathbf{U}, x_{0}, y_{0})$, whose arguments denote a solution $\mathbf{U}$ to the shallow water equations and the $x$- and $y$-positions of a seed point, respectively.
	A complete description of this function is given in the appendix (Algorithm \ref{alg:streamline}).

	The definition of $\textproc{ComputePathline}$ enables the computation of a set of points $Q$ approximating the \emph{pathtube} leading to a set of exterior asset points $P_{a} \in a \in \mathcal{A}$ via
	\begin{linenomath*}
	\begin{equation}
	Q(\mathbf{U}, P_{a}) = \bigcup_{\mathclap{(x_{0}, y_{0}) \in P_{a}}} \textproc{ComputePathline}(\mathbf{U}, x_{0}, y_{0}).
	\end{equation}
	\end{linenomath*}
	Since pathtubes are curvilinear, typical geometries that envelope $Q$ (e.g., the convex hull) do not effectively summarize this set.
	For this reason, the notion of an \emph{alpha shape} is introduced, which minimally encompasses points of $Q$ using straight lines.
	A discussion on alpha shapes can be found in \citeA{alphashapesintro}.
	In this study, Edelsbrunner's algorithm (\citeA{edelsbrunner1983shape}), presented in the appendix (Algorithm \ref{alg:edelsbrunner}), is used to compute alpha shapes.
	The function that computes this shape for a set $Q$ and alpha value $\alpha$ is denoted as $\textproc{AlphaShape}(Q, \alpha)$.

	The definition of the function $\textproc{AlphaShape}$ finally allows for definition of the functions $\textproc{InitializeRestriction}$ and $\textproc{UpdateRestriction}$.
	Both assume restrictions are the unions of alpha shapes approximating the pathtubes leading to each asset.
	The functions are described in Algorithms \ref{alg:initialize-restricted-region-streamlines} and Equation \eqref{eqn:update-restriction}, respectively.
	In Algorithm \ref{alg:initialize-restricted-region-streamlines}, Line \ref{alg:initialize-restricted-region-streamlines:line:initialize}, the shallow water equations are solved \emph{without} the presence of structural mitigation measures.
	In Line \ref{alg:initialize-restricted-region-streamlines:line:get-p}, the union of alpha shapes for all pathtubes leading to the assets $a \in \mathcal{A}$ is computed.
	Asset regions are then subtracted from this set to ensure structural measures do not overlap with asset locations.

	The function \textproc{UpdateRestriction} using the pathline approach is defined as
	\begin{linenomath*}
	\begin{equation}
	\textproc{UpdateRestriction}(\mathbf{U}, \mathcal{A}, \mathcal{P}, \alpha) = \mathcal{P} \cup \left(\bigcup_{a \in \mathcal{A}} \textproc{AlphaShape}(Q(\mathbf{U}, P_{a}), \alpha)\right) \setminus \bigcup \mathcal{A} \label{eqn:update-restriction}.
	\end{equation}
	\end{linenomath*}
	The majority of this function resembles Algorithm \ref{alg:initialize-restricted-region-streamlines}, although the union of the current set and previous $\mathcal{P}$ is computed to encourage exploration of a more representative (i.e., expanded) search space.
	Moreover, as per Algorithm \ref{alg:pseudoalg}, this function is only called as better solutions to the OFMP are obtained.
	This decreases the burden of computing pathtubes and alpha shapes on each iteration of the algorithm.

	\begin{algorithm}[t]
		 \caption{\textproc{InitializeRestriction}: Returns the initial restricted positional set.}
		 \label{alg:initialize-restricted-region-streamlines}
		 \begin{algorithmic}[1]
			  \Function{InitializeRestriction}{$B, n, \mathcal{A}, \alpha$}
					\State $\textrm{Solve} ~ \mathbf{U}_{t} + \mathbf{F}_{x} + \mathbf{G}_{y} = \mathbf{S}_{R} + \mathbf{S}_{B} + \mathbf{S}_{n}$ \label{alg:initialize-restricted-region-streamlines:line:initialize}
					\State \Return $\bigcup_{a \in \mathcal{A}} \textproc{AlphaShape}\left(Q(\mathbf{U}, P_{a}), \alpha\right) \setminus \bigcup \mathcal{A}$ \label{alg:initialize-restricted-region-streamlines:line:get-p}
			  \EndFunction
		 \end{algorithmic}
	\end{algorithm}

	\subsection{Sequential optimization algorithm}
	\label{sec:algorithm-hierarchical}
	Due to the nonlinear sensitivity of flooding behavior with respect to mitigation efforts, predictable and incremental changes to solutions of the OFMP while increasing the number of mitigation measures, $m$, are not ensured.
	This may be undesirable from a planning perspective.
	A separate algorithm is thus proposed to induce a sequential solution to the OFMP, whereby solutions with $m = 2$ include those of $m = 1$, solutions with $m = 3$ include those of $m = 2$, and so on.
	This ensures increasing utility for configurations of increasing sizes.
	It also allows policymakers to more clearly understand the effects of budgetary constraints with respect to the overall structural flood mitigation efforts.
	The recursion to compute sequential solutions may be defined as
	\begin{linenomath*}
	\begin{equation}
	(B_{i}, n_{i}) = \textproc{SolveOFMP}\left(B_{i - 1}, n_{i - 1}, \mathcal{A}, 1, \frac{T_{\max}}{m}, \alpha\right),
	\end{equation}
	\end{linenomath*}
	where $B_{0} = B$, $n_{0} = n$, and the time for each subproblem is an equal portion of $T_{\max}$.
	In this paper, Line \ref{alg:pseudoalg:line:update-restricted-region} is eliminated from Algorithm \ref{alg:pseudoalg} when using the sequential approach, as the best placement for a \emph{single} structural mitigation measure is likely to reside within the \emph{initial} $\tilde{\mathcal{P}}$ computed on Line \ref{alg:pseudoalg:line:init-2}.
	As a consequence, for each structural measure placed using the sequential approach, pathtubes are constructed only once.

	\section{Results}
	\label{sec:results}

	\subsection{Model relaxation}
	\label{subsec:experimental_setting}
	The proposed approach uses the open source \path{scipy.optimize.differential_evolution} (DE) and RBFOpt libraries to produce two separate implementations of \textproc{GenerateSolution} in Algorithm \ref{alg:pseudoalg} \cite{storn1997differential, costa2014rbfopt}.
	Both only include support for simple bounds like those indicated in Constraints \eqref{eqn:specific-ofmp-valid-structures}.
	Thus, these implementations of \textproc{GenerateSolution} may generate configurations that are infeasible with respect to Constraints \eqref{eqn:specific-ofmp-no-overlap} through \eqref{eqn:specific-ofmp-wall-definition-veg}.
	To overcome this, the OFMP defined by Equations \eqref{eqn:specific-ofmp-objective} through \eqref{eqn:specific-ofmp-centroid-restriction} is replaced with the relaxed formulation
	\begin{linenomath*}
	\begin{subequations}
	\begin{align}
		 & \underset{\omega_{1}, \omega_{2}, \dots, \omega_{m}}{\text{minimize}}
		 & & z\left(\omega_{1}, \omega_{2}, \dots, \omega_{m}\right) = p_{1} + p_{2} + \mathlarger{\sum}_{a \in \mathcal{A}} \iint_{a} \max_{t} {h(x, y, t)} \,dx \,dy \label{eqn:specific-ofmp-objective-mod} \\
		 & \text{subject to}
		 & & p_{1} = c_{1} \mathlarger{\sum}_{i = 1}^{m} \min_{} \left\{\left\lVert(x, y) - \left(\lambda_{i}, \phi_{i}\right)\right\rVert : (x, y) \in \mathcal{P}\right\} \label{eqn:specific-ofmp-penalty-mod} \\
		 & & & p_{2} = c_{2} \mathlarger{\sum}_{i = 1}^{m} \mathlarger{\sum}_{a \in \mathcal{A}} \iint_{a} \delta_{B}(\omega_{i}) \,dx \,dy + c_{3} \mathlarger{\sum}_{i = 1}^{m} \mathlarger{\sum}_{a \in \mathcal{A}} \iint_{a} \delta_{n}(\omega_{i}) \,dx \,dy \label{eqn:specific-ofmp-penalty-mod-1} \\
		 & & & \mathbf{U}_{t} + \mathbf{F}_{x} + \mathbf{G}_{y} = \mathbf{S}_{R} + \tilde{\mathbf{S}}_{B} + \tilde{\mathbf{S}}_{n} \label{eqn:specific-ofmp-2d-swe-mod} \\
		 & & & \delta_{B}(\omega_{i})(x, y) =
			  \begin{cases}
			  \begin{aligned}
				  \bar{b}_{i} & ~ \textrm{for}
						 \begin{cases}
							  |(x - \phi_{i}) \cos \theta_{i} - (y - \lambda_{i}) \sin \theta_{i}| \leq \frac{\ell}{2} \\
							  |(x - \phi_{i}) \sin \theta_{i} + (y - \lambda_{i}) \cos \theta_{i}| \leq \frac{w}{2}
						 \end{cases} \\
			  0 & ~ \textrm{otherwise}
			  \end{aligned}
		  \end{cases} \hspace{-2.05em} \forall i \in \mathcal{M} \label{eqn:specific-ofmp-wall-definition-mod} \\
		 & & & \delta_{n}(\omega_{i})(x, y) =
			  \begin{cases}
				  \bar{n}_{i} & \textrm{for} ~ (x - \phi_{i})^{2} + (y - \lambda_{i})^{2} \leq r^{2} \\
				  0 & \textrm{otherwise}
			  \end{cases} \forall i \in \mathcal{M} \label{eqn:specific-ofmp-wall-definition-veg-mod} \\
		 & & & \lambda_{lb} \leq \lambda_{i} \leq \lambda_{ub}, ~ \phi_{lb} \leq \phi_{i} \leq \phi_{ub}, ~ 0 \leq \theta_{i} \leq \pi, ~ \forall i \in \mathcal{M} \label{eqn:specific-ofmp-valid-structures-mod}.
	\end{align}
	\end{subequations}
	\end{linenomath*}
	In Equation \eqref{eqn:specific-ofmp-objective-mod}, two penalty terms are included in the objective to capture infeasibilities in Constraints \eqref{eqn:specific-ofmp-no-overlap} through \eqref{eqn:specific-ofmp-wall-definition-veg}.
	The first penalty, $p_{1}$, is defined in Constraint \eqref{eqn:specific-ofmp-penalty-mod} and denotes the sum of all minimum distances between each measure's centroid and the nearest point of the restricted positional set $\mathcal{P}$.
	This term is scaled by the constant $c_{1}$, taken in this study to be equal to one.
	The second penalty, $p_{2}$, is defined in Constraint \eqref{eqn:specific-ofmp-penalty-mod-1}.
	Here, the first term denotes the net modified elevation volume over all asset regions, and the second term denotes the net change in roughness over all asset regions.
	These terms are scaled by the constants $c_{2}$ and $c_{3}$, respectively.
	Herein, both are taken to be $(\Delta r)^{-2}$, where $\Delta r$ is the spatial resolution of the discretization.

	\subsection{Experimental setting}
	\label{sec:experimental-setting}
	For simplicity, Sections \ref{subsect:results-simple} through \ref{sec:alg-summary} focus on OFMPs where only wall-type measures are considered (i.e., $\bar{b}_{i} > 0$), while Section \ref{sec:proof-of-concept} presents an algorithmic proof of concept where only revegetation-type measures are considered (i.e., $\bar{n}_{i} > 0$).
	For each experiment, Algorithm \ref{alg:pseudoalg} was limited to one day of wall-clock time.
	When using DE, population sizes of $45 m$ ($\bar{b}_{i} > 0$) and $30 m$ ($\bar{n}_{i} > 0$) were employed; trial solutions were computed as the best solution plus scaled contributions of two random candidates; the mutation constant varied randomly within $[0.5, 1.0)$; and the recombination constant was set to $0.9$.
	When using the direct $\textproc{InitializeRestriction}$ and $\textproc{UpdateRestriction}$ methods, Latin hypercube sampling was used to initialize the population.
	When using the pathline-based methods, the population was initialized via random sampling over the initial restricted set (i.e., $\tilde{\mathcal{P}}$) and $\theta_{i} \in [0, \pi]$.
	When using RBFOpt, the \texttt{sampling} method was used; most other parameters were left unchanged.

	For computational considerations, if the configuration proposed by \textproc{GenerateSolution} was feasible, the shallow water equations (i.e., Constraint \eqref{eqn:specific-ofmp-2d-swe-mod}) were solved using the proposed configuration.
	Otherwise, a solution containing \emph{no structural mitigation measures} was referenced.
	That is, $\tilde{S}_{B}$ was replaced with $S_{B}$, and $\tilde{S}_{n}$ was replaced with $S_{n}$.
	To solve these PDEs, the open-source surface-water modeling software Nuflood \cite{tasseff2016nuflood} was used, where the shallow water equations are spatially discretized according to the scheme described by \citeA{kurganov2007second}.

	Each experiment was conducted on one Intel Xeon E5-2695 V4 CPU containing eighteen cores at 2.1 GHz and 125 GB of RAM.
	Nuflood was compiled in single-precision mode using the Intel C\texttt{++} Compiler, version 17.0.1.
	The remainder of Algorithm \ref{alg:pseudoalg} was implemented in Python 3.6.
	Compared to the PDE evaluations, these other portions of the algorithm were found to be computationally negligible.

	\subsection{Simplified circular dam break scenarios}
	\label{subsect:results-simple}
	\begin{figure}[t]
		\includegraphics[width=\textwidth]{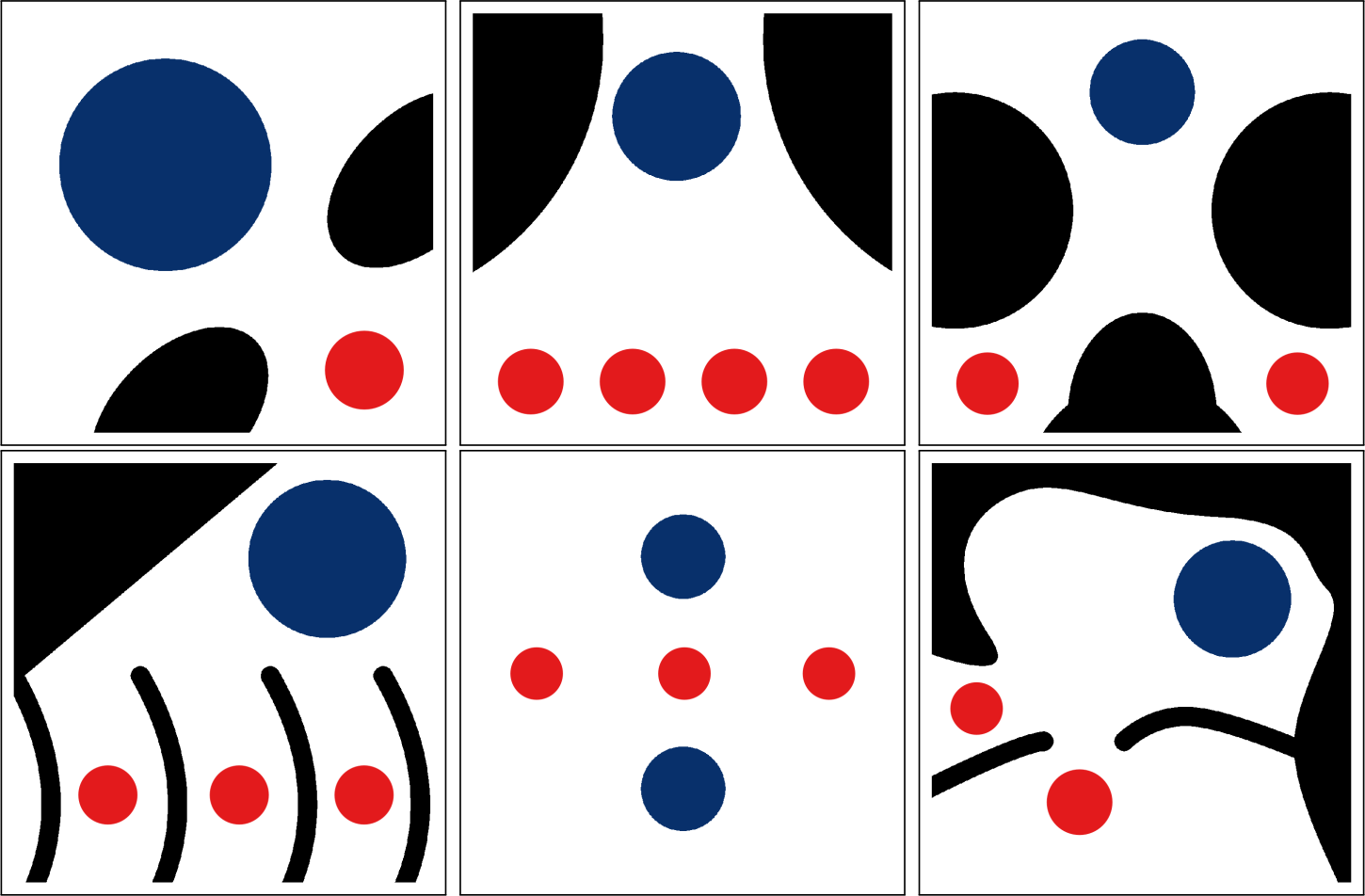}
		\caption{Pictorial descriptions of six simple OFMP scenarios, ordered numerically (e.g., one in the upper left). Black represents nonzero topographic elevation (of height one meter); blue represents nonzero initial water depth (of height one meter); and red represents assets.}
		\label{fig:simple-scenarios}
		\vspace{-3em}
	\end{figure}

	\begin{figure}[!ht]
		\includegraphics[width=\textwidth]{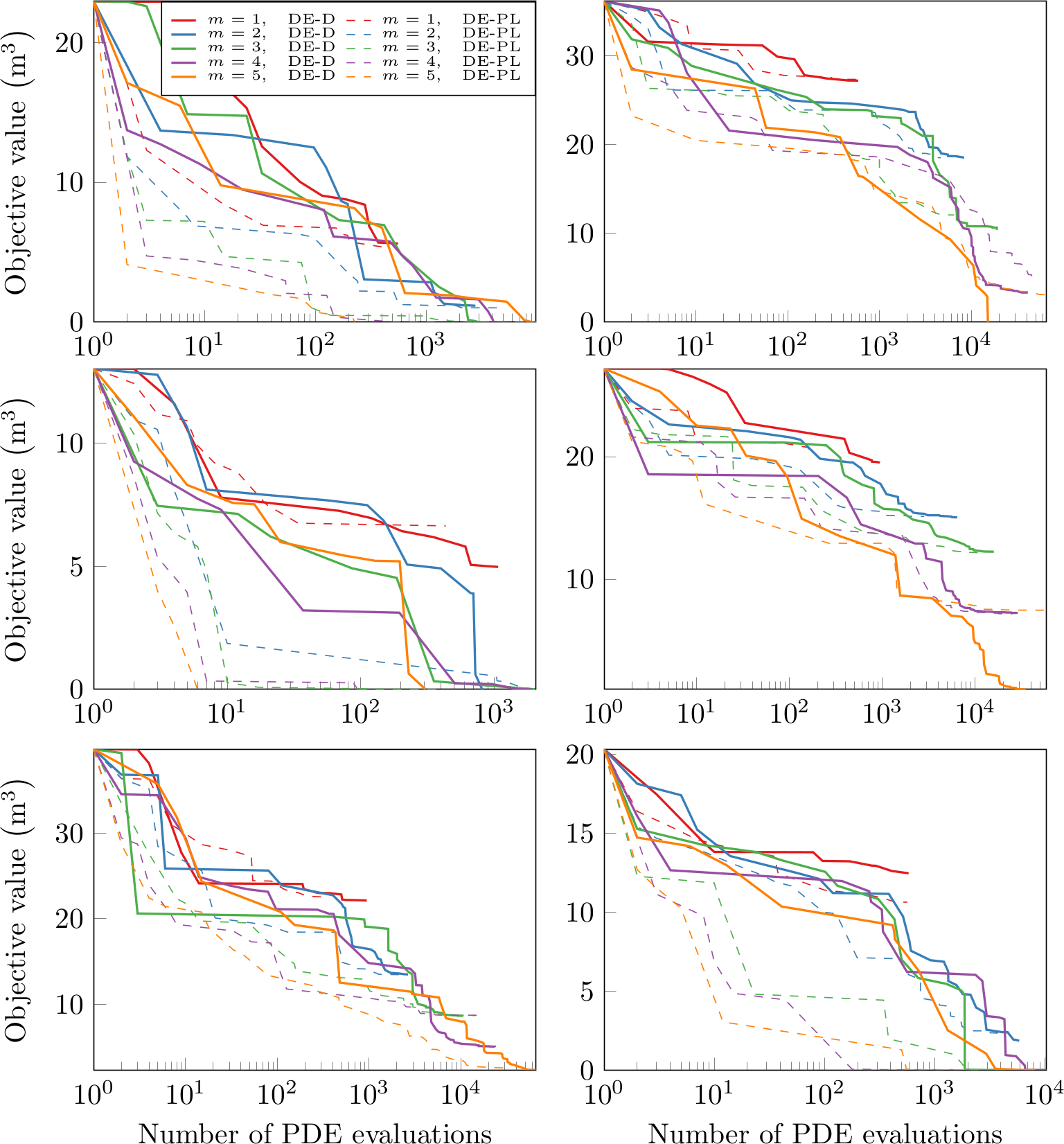}
		\caption{Comparison of objective value versus number of PDE evaluations for Scenarios 1 through 6, respectively, using DE-D and DE-PL for configurations of one through five walls.}
		\label{fig:simple-objective-plots}
		\vspace{-3em}
	\end{figure}

	\begin{figure}[!ht]
		\includegraphics[width=\textwidth]{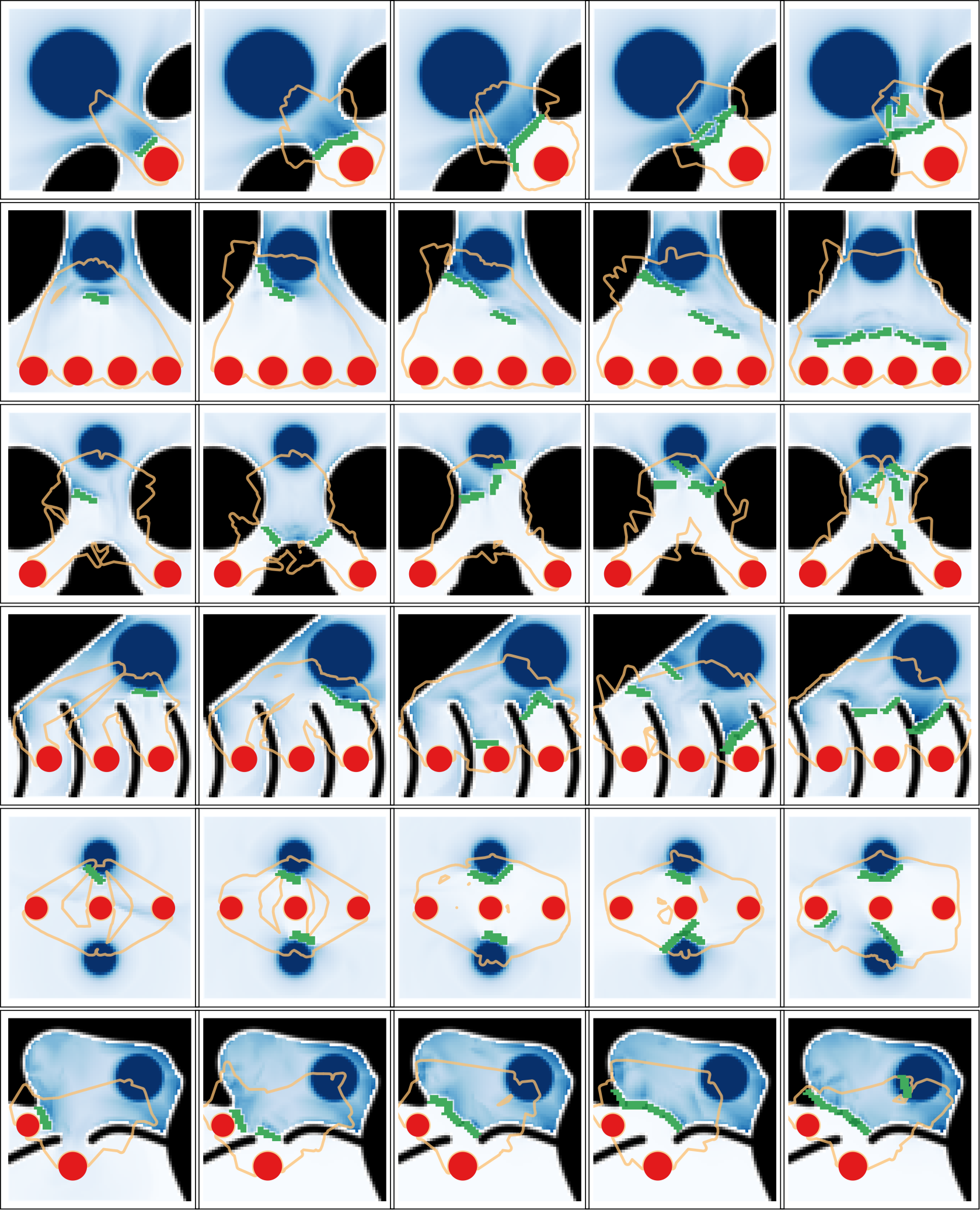}
		\caption{Best obtained elevations and maximum depths for configurations of one through five walls for the highly simplified flood scenarios, referenced as Scenarios 1 through 6 (vertically). Darker blue corresponds to larger maximum depths; black corresponds to nonzero portions of the initial topographic elevation field; green corresponds to elevation additions via the placement of walls; and red corresponds to asset locations. The orange lines represent the exteriors of the final computed restricted positional sets $\tilde{\mathcal{P}}$ in Algorithm \ref{alg:pseudoalg}.}
		\label{fig:best-simple-summaries}
		\vspace{-3em}
	\end{figure}

	To compare the two positional restriction methodologies described in Section \ref{sec:algorithm-initialize-update-restricted-region}, six simple OFMP scenarios were constructed.
	All were intended to have human intuitive solutions, i.e., optimal placement of structural mitigation measures could be inferred from a basic understanding of flood propagation.
	These scenarios are displayed pictorially in Figure \ref{fig:simple-scenarios}.
	In each scenario, under the influence of gravity, the initial volume of water (colored with blue) is propagated outward; without mitigation measures, this water comes into contact with assets (colored with red), flooding them.

	Each of the six scenarios was modeled using a spatial resolution of one meter and $64 \times 64$ grid cells.
	The ground surface was assumed to be frictionless; critical depth boundary conditions were employed; and a simulation duration of one hundred seconds was used.
	When necessary to compute pathlines, intermediate PDE solution data was reported for every one second of simulation time.
	In the experiments performed, each of the corresponding OFMPs was solved with the number of walls, $m$, ranging from one to five.
	Wall widths, lengths, and heights were fixed to 2.5, 8.0, and 1.0 meters, respectively.
	Finally, all experiments were performed using a single fixed random seed.

	In Figure \ref{fig:simple-objective-plots}, for each experiment, the objective behavior is plotted against the number of PDE evaluations required to reach that objective.
	These behaviors are compared for the direct differential evolution solver (DE-D) and its pathline-based counterpart (DE-PL).
	The DE-PL solver was generally able to find good solutions faster and improve upon them more rapidly, especially for configurations involving larger numbers of walls.
	However, there were some instances where the DE-D solver produced higher quality solutions than the DE-PL solver, e.g., when optimizing the configuration of five walls in Scenario 4.
	These anomalies could be a consequence of the random nature of the DE algorithm; they could also be due to the DE implementation's tendency to terminate once a population has sufficiently stabilized.

	In Figure \ref{fig:best-simple-summaries}, the best obtained wall configurations using DE-PL are displayed pictorially for all pairs of scenarios and numbers of mitigation measures.
	The configurations resemble what might be intuited by a human.
	When applicable, configurations are non-overlapping and well-connected.
	As the number of walls varies, configurations also show interesting nonincremental behavior.
	For example, in the first scenario, walls are initially placed close to the asset; as the number of walls increases, they are placed farther away to form connections with existing topographic features.
	However, such non-sequential behavior may be undesirable from a planning perspective.

	\subsection{Hypothetical dam break scenario from Theme C of the 12th International Benchmark Workshop on Numerical Analysis of Dams}
	\label{subsect:results-icold}
	This section focuses on demonstrating the merits of the sequential optimization algorithm using the hypothetical dam break defined in Theme C of the 12th International Benchmark Workshop on Numerical Analysis of Dams (ICOLD 2013) \cite{judi2014computational}.
	To simulate this scenario, the dam break was modeled as a point source with time-dependent discharge.
	The initial topographic elevation field (with the dam excluded) was provided by the workshop and resampled from a resolution of ten to ninety meters to ease computational burden.
	The Manning's roughness coefficient was set to $0.035$; critical depth boundary conditions were employed; a duration of twelve hours was used; and, when necessary to compute pathlines, PDE solution data was reported every ten minutes of simulation time.

	Asset locations and sizes were selected to increase the difficulty of the OFMP, with two assets placed near the primary channel of the scenario and three placed farther away. The experimental setup remained similar to that described in Section \ref{subsect:results-simple}.
	However, in this case, the number of walls ranged from one to ten, while wall widths, lengths, and heights were fixed to 250, 1000, and ten meters, respectively.
	To compare differences in OFMP solver performance, each solver was executed using ten different random seeds for each possible number of walls.
	In total, the experiments described in this subsection thus required nearly six hundred days of compute time.

	\subsubsection{Pathline-based algorithm results}
	\begin{table}[t]
		 \begin{tabular}{c|c|c|c|c|c|c|c|c|c|}
			\cline{2-10} & \multicolumn{4}{c|}{RBFOpt-D} & \multicolumn{4}{c|}{RBFOpt-PL} & \multirowcell{2}{Mean \\ Improvement} \\
			\cline{1-9} \multicolumn{1}{|c|}{$m$} & Mean & Min & Max & SD & Mean & Min & Max & SD & \\
			\cline{1-10} \multicolumn{1}{|c|}{$1$} & $165.93$ & $159.63$ & $170.25$ & $3.36$ & $162.81$ & $159.63$ & $167.59$ & $2.61$ & $1.88\%$ \\
			\cline{1-10} \multicolumn{1}{|c|}{$2$} & $147.96$ & $105.34$ & $166.60$ & $20.14$ & $111.02$ & $95.98$ & $130.88$ & $10.74$ & $24.96\%$ \\
			\cline{1-10} \multicolumn{1}{|c|}{$3$} & $144.93$ & $111.87$ & $164.10$ & $15.44$ & $89.22$ & $77.56$ & $93.62$ & $5.40$ & $38.44\%$ \\
			\cline{1-10} \multicolumn{1}{|c|}{$4$} & $128.52$ & $105.81$ & $159.12$ & $14.21$ & $81.33$ & $51.63$ & $97.83$ & $13.47$ & $36.72\%$ \\
			\cline{1-10} \multicolumn{1}{|c|}{$5$} & $135.94$ & $128.25$ & $144.42$ & $5.03$ & $62.38$ & $26.40$ & $80.40$ & $15.67$ & $54.11\%$ \\
			\cline{1-10} \multicolumn{1}{|c|}{$6$} & $122.24$ & $98.38$ & $140.19$ & $14.22$ & $59.13$ & $41.40$ & $71.25$ & $9.07$ & $51.63\%$ \\
			\cline{1-10} \multicolumn{1}{|c|}{$7$} & $119.14$ & $81.14$ & $145.80$ & $18.65$ & $51.11$ & $29.13$ & $66.20$ & $13.08$ & $57.10\%$ \\
			\cline{1-10} \multicolumn{1}{|c|}{$8$} & $102.08$ & $78.65$ & $122.69$ & $16.22$ & $43.28$ & $27.82$ & $57.35$ & $9.40$ & $57.60\%$ \\
			\cline{1-10} \multicolumn{1}{|c|}{$9$} & $107.40$ & $81.93$ & $127.53$ & $15.52$ & $42.73$ & $19.00$ & $53.26$ & $11.19$ & $60.21\%$ \\
			\cline{1-10} \multicolumn{1}{|c|}{$10$} & $104.90$ & $77.17$ & $124.13$ & $16.75$ & $34.72$ & $21.75$ & $42.74$ & $7.60$ & $66.90\%$ \\
			\cline{1-10}
		 \end{tabular}
		 \vspace{2em}
		 \caption{Table comparing objectives using the (direct) RBFOpt-D and (pathline-based) RBFOpt-PL solvers over ten random seeds, with the number of walls ($m$) ranging from one to ten, as discussed in Section \ref{subsect:results-icold}. Values are scaled by a factor of $10^{-4}$.}
		 \label{tab:icold-stats-table-1}
	\end{table}

	\begin{table}[t]
		 \begin{tabular}{c|c|c|c|c|c|c|c|c|c|}
			\cline{2-10} & \multicolumn{4}{c|}{DE-D} & \multicolumn{4}{c|}{DE-PL} & \multirowcell{2}{Mean \\ Improvement} \\
			\cline{1-9} \multicolumn{1}{|c|}{$m$} & Mean & Min & Max & SD & Mean & Min & Max & SD & \\
			\cline{1-10} \multicolumn{1}{|c|}{$1$} & $\underline{162.18}$ & $\mathbf{158.59}$ & $167.59$ & $3.85$ & $163.86$ & $159.63$ & $170.07$ & $4.52$ & $-1.04\%$ \\
			\cline{1-10} \multicolumn{1}{|c|}{$2$} & $99.49$ & $\mathbf{84.87}$ & $104.24$ & $6.29$ & $102.61$ & $94.24$ & $105.39$ & $3.19$ & $-3.13\%$ \\
			\cline{1-10} \multicolumn{1}{|c|}{$3$} & $66.78$ & $39.56$ & $119.20$ & $27.76$ & $\underline{50.75}$ & $\mathbf{33.36}$ & $65.93$ & $14.80$ & $24.01\%$ \\
			\cline{1-10} \multicolumn{1}{|c|}{$4$} & $101.06$ & $79.57$ & $134.61$ & $20.57$ & $31.84$ & $17.51$ & $56.87$ & $13.20$ & $68.49\%$ \\
			\cline{1-10} \multicolumn{1}{|c|}{$5$} & $115.43$ & $100.56$ & $145.39$ & $16.04$ & $22.74$ & $12.82$ & $36.39$ & $8.04$ & $80.30\%$ \\
			\cline{1-10} \multicolumn{1}{|c|}{$6$} & $124.15$ & $101.22$ & $145.75$ & $15.85$ & $24.79$ & $7.24$ & $60.66$ & $14.87$ & $80.03\%$ \\
			\cline{1-10} \multicolumn{1}{|c|}{$7$} & $129.17$ & $110.94$ & $153.14$ & $12.52$ & $18.14$ & $5.59$ & $26.77$ & $6.64$ & $85.96\%$ \\
			\cline{1-10} \multicolumn{1}{|c|}{$8$} & $125.51$ & $96.68$ & $144.92$ & $16.47$ & $14.02$ & $4.03$ & $19.74$ & $5.59$ & $88.83\%$ \\
			\cline{1-10} \multicolumn{1}{|c|}{$9$} & $129.23$ & $99.39$ & $146.58$ & $14.02$ & $17.83$ & $8.56$ & $22.46$ & $4.85$ & $86.20\%$ \\
			\cline{1-10} \multicolumn{1}{|c|}{$10$} & $119.09$ & $86.16$ & $140.87$ & $17.91$ & $19.68$ & $10.43$ & $30.38$ & $5.61$ & $83.47\%$ \\
			\cline{1-10}
		 \end{tabular}
		 \vspace{2em}
		 \caption{Table comparing objective values obtained using the (direct) DE-D and (pathline-based) DE-PL solvers over ten random seeds, with the number of walls ($m$) ranging from one to ten. Values are scaled by a factor of $10^{-4}$. Best objectives over all seeds and solvers in Tables \ref{tab:icold-stats-table-1} through \ref{tab:icold-stats-table-4} are denoted in bold, while best \emph{mean} objectives are underlined.}
		 \label{tab:icold-stats-table-2}
		 \vspace{-3em}
	\end{table}

	To confirm the effectiveness of the pathline-based solvers, two implementations of Algorithm \ref{alg:pseudoalg} using RBFOpt were benchmarked.
	In Table \ref{tab:icold-stats-table-1}, the objective behavior of the pathline-based solver (RBFOpt-PL) is compared against its direct counterpart (RBFOpt-D).
	The pathline-based solver clearly outperformed RBFOpt-D in nearly all instances, e.g., it resulted in smaller minima, means, and standard deviations.
	The single exception appears to be for $m = 1$, where the direct solver produced an equivalent minimum to the pathline-based solver.
	Nonetheless, on average, the pathline-based solver provided a $45\%$ improvement over the direct solver, with generally larger improvements for greater numbers of walls.
	This improvement was computed as
	\begin{linenomath*}
	\begin{equation}
		\textrm{percentage improvement} = 100 \left(\frac{a - b}{a}\right),
	\end{equation}
	\end{linenomath*}
	where, here, $a$ and $b$ represent the mean objective values obtained from the RBFOpt-D and RBFOpt-PL solvers.
	The same metric is also used throughout Tables \ref{tab:icold-stats-table-2}, \ref{tab:icold-stats-table-3}, and \ref{tab:icold-stats-table-4}.

	A similar comparison is made between DE-D and DE-PL in Table \ref{tab:icold-stats-table-2}.
	Again, the pathline-based solver (DE-PL) outperformed its direct counterpart (DE-D) in nearly all metrics, providing an overall mean improvement of $59\%$.
	The pathline-based solver also displayed mostly monotonic decreases in the objective as the number of walls increased, while the objectives associated with the direct solver generally \emph{increased} as the number of walls increased.
	However, note that for small numbers of walls (i.e., one and two), the direct DE solver outperformed its pathline-based counterpart.
	This could be a consequence of the more complicated objective penalty in Constraint \eqref{eqn:specific-ofmp-penalty-mod} when $\mathcal{P}$ is restricted.
	For example, \citeA{DEB2000311} describes various means by which penalty-based genetic algorithms can result in nonoptimal solutions.
	Nonetheless, overall, the direct penalization method considered herein works well.

	It is important to note the differences between the RBFOpt-based and DE-based solvers benchmarked in Tables \ref{tab:icold-stats-table-1} and \ref{tab:icold-stats-table-2}, respectively.
	In general, DE-PL greatly outperformed both RBFOpt-based solvers; for example, DE-PL provided a $47\%$ mean improvement over RBFOpt-PL.
	These differences could be for multiple reasons.
	For example, there are many more hyperparameters associated with RBFOpt than DE; more careful tuning may have increased RBFOpt's convergence.
	Furthermore, RBFOpt's \texttt{sampling} search strategy was used to show the efficacy of the pathline-based approach when applied to other (non-evolutionary) search techniques; the solver software may have performed more favorably using some other strategy.

	Figure \ref{fig:best-icold-summaries} displays the best obtained wall configuration for each possible number of walls using the DE-PL solver.
	Structure placement appears highly nonincremental as the number of walls increases, especially for smaller numbers of walls.
	Also, when optimizing for a number of walls greater than eight, solutions generally deteriorated, indicating the search space becomes prohibitively large.
	Interestingly, the size of the restricted set $\tilde{\mathcal{P}}$ did not increase substantially as the configuration size grew.
	Finally, in Figure \ref{fig:de-icold-10}, the best obtained solution for ten walls using DE-D is displayed; this underscores the difficulty of such a problem when applying a conventional algorithm.

	\subsubsection{Sequential algorithm results}
	\begin{table}[t]
		 \begin{tabular}{c|c|c|c|c|c|c|c|c|c|}
			  \cline{2-10} & \multicolumn{4}{c|}{DE-PL} & \multicolumn{4}{c|}{DE-D-S} & \multirowcell{2}{Mean \\ Improvement} \\
			  \cline{1-9} \multicolumn{1}{|c|}{$m$} & Mean & Min & Max & SD & Mean & Min & Max & SD & \\
			\cline{1-10} \multicolumn{1}{|c|}{$1$} & $163.86$ & $159.63$ & $170.07$ & $4.52$ & $\underline{162.18}$ & $\mathbf{158.59}$ & $167.59$ & $3.85$ & $1.02\%$ \\
			\cline{1-10} \multicolumn{1}{|c|}{$2$} & $102.61$ & $94.24$ & $105.39$ & $3.19$ & $102.38$ & $100.97$ & $103.68$ & $0.87$ & $0.22\%$ \\
			\cline{1-10} \multicolumn{1}{|c|}{$3$} & $\underline{50.75}$ & $\mathbf{33.36}$ & $65.93$ & $14.80$ & $64.98$ & $48.97$ & $81.93$ & $14.71$ & $-28.06\%$ \\
			\cline{1-10} \multicolumn{1}{|c|}{$4$} & $31.84$ & $17.51$ & $56.87$ & $13.20$ & $41.04$ & $26.07$ & $58.60$ & $14.88$ & $-28.89\%$ \\
			\cline{1-10} \multicolumn{1}{|c|}{$5$} & $22.74$ & $12.82$ & $36.39$ & $8.04$ & $23.22$ & $17.45$ & $34.87$ & $6.18$ & $-2.11\%$ \\
			\cline{1-10} \multicolumn{1}{|c|}{$6$} & $24.79$ & $7.24$ & $60.66$ & $14.87$ & $14.94$ & $11.56$ & $18.39$ & $2.26$ & $39.76\%$ \\
			\cline{1-10} \multicolumn{1}{|c|}{$7$} & $18.14$ & $5.59$ & $26.77$ & $6.64$ & $11.14$ & $4.10$ & $14.53$ & $2.87$ & $38.60\%$ \\
			\cline{1-10} \multicolumn{1}{|c|}{$8$} & $14.02$ & $4.03$ & $19.74$ & $5.59$ & $9.42$ & $6.02$ & $16.64$ & $3.30$ & $32.81\%$ \\
			\cline{1-10} \multicolumn{1}{|c|}{$9$} & $17.83$ & $8.56$ & $22.46$ & $4.85$ & $4.72$ & $\mathbf{0.00}$ & $9.91$ & $3.74$ & $73.52\%$ \\
			\cline{1-10} \multicolumn{1}{|c|}{$10$} & $19.68$ & $10.43$ & $30.38$ & $5.61$ & $3.07$ & $\mathbf{0.00}$ & $9.29$ & $3.10$ & $84.38\%$ \\
			  \cline{1-10}
		 \end{tabular}
		 \vspace{2em}
		 \caption{Table comparing objective values obtained using the (pathline-based) DE-PL and (direct sequential) DE-D-S solvers over ten random seeds, with the number of walls ($m$) ranging from one to ten. Values are scaled by a factor of $10^{-4}$. Best objectives over all seeds and solvers in Tables \ref{tab:icold-stats-table-1} through \ref{tab:icold-stats-table-4} are denoted in bold, while best \emph{mean} objectives are underlined.}
		 \label{tab:icold-stats-table-3}
	\end{table}

	\begin{table}[t]
		 \begin{tabular}{c|c|c|c|c|c|c|c|c|c|}
			\cline{2-10} & \multicolumn{4}{c|}{DE-D-S} & \multicolumn{4}{c|}{DE-PL-S} & \multirowcell{2}{Mean \\ Improvement} \\
			\cline{1-9} \multicolumn{1}{|c|}{$m$} & Mean & Min & Max & SD & Mean & Min & Max & SD & \\
			\cline{1-10} \multicolumn{1}{|c|}{$1$} & $\underline{162.18}$ & $\mathbf{158.59}$ & $167.59$ & $3.85$ & $163.61$ & $159.63$ & $167.59$ & $4.19$ & $-0.88\%$ \\
			\cline{1-10} \multicolumn{1}{|c|}{$2$} & $102.38$ & $100.97$ & $103.68$ & $0.87$ & $\underline{98.25}$ & $86.78$ & $105.84$ & $7.45$ & $4.04\%$ \\
			\cline{1-10} \multicolumn{1}{|c|}{$3$} & $64.98$ & $48.97$ & $81.93$ & $14.71$ & $56.74$ & $35.73$ & $86.12$ & $13.92$ & $12.69\%$ \\
			\cline{1-10} \multicolumn{1}{|c|}{$4$} & $41.04$ & $26.07$ & $58.60$ & $14.88$ & $\underline{27.18}$ & $\mathbf{14.55}$ & $56.37$ & $11.81$ & $33.77\%$ \\
			\cline{1-10} \multicolumn{1}{|c|}{$5$} & $23.22$ & $17.45$ & $34.87$ & $6.18$ & $\underline{16.32}$ & $\mathbf{8.38}$ & $25.35$ & $5.25$ & $29.73\%$ \\
			\cline{1-10} \multicolumn{1}{|c|}{$6$} & $14.94$ & $11.56$ & $18.39$ & $2.26$ & $\underline{10.63}$ & $\mathbf{3.92}$ & $16.96$ & $5.29$ & $28.82\%$ \\
			\cline{1-10} \multicolumn{1}{|c|}{$7$} & $11.14$ & $4.10$ & $14.53$ & $2.87$ & $\underline{7.02}$ & $\mathbf{0.13}$ & $15.00$ & $5.08$ & $36.97\%$ \\
			\cline{1-10} \multicolumn{1}{|c|}{$8$} & $9.42$ & $6.02$ & $16.64$ & $3.30$ & $\underline{4.53}$ & $\mathbf{0.00}$ & $9.68$ & $3.81$ & $51.90\%$ \\
			\cline{1-10} \multicolumn{1}{|c|}{$9$} & $4.72$ & $\mathbf{0.00}$ & $9.91$ & $3.74$ & $\underline{3.36}$ & $\mathbf{0.00}$ & $7.86$ & $3.43$ & $28.83\%$ \\
			\cline{1-10} \multicolumn{1}{|c|}{$10$} & $3.07$ & $\mathbf{0.00}$ & $9.29$ & $3.10$ & $\underline{2.59}$ & $0.00$ & $6.30$ & $2.60$ & $15.76\%$ \\
		  \cline{1-10}
		 \end{tabular}
		 \vspace{2em}
		 \caption{Table comparing objectives obtained using the (direct sequential) DE-D-S and (pathline-based sequential) DE-PL-S solvers over ten random seeds, with the number of walls ($m$) ranging from one to ten. Values are scaled by a factor of $10^{-4}$. Best objectives over all seeds and solvers in Tables \ref{tab:icold-stats-table-1} through \ref{tab:icold-stats-table-4} are denoted in bold, while best \emph{mean} objectives are underlined.}
		 \label{tab:icold-stats-table-4}
		 \vspace{-3em}
	\end{table}

	To counteract the degradation of solutions for larger configurations, the sequential approach presented in Section \ref{sec:algorithm-hierarchical} was benchmarked in a similar setting.
	In Table \ref{tab:icold-stats-table-3}, performance of the direct sequential DE solver (DE-D-S) is compared against DE-PL.
	Interestingly, DE-D-S performed much better than DE-PL for configurations containing many walls, providing improvements as large as $84\%$. This result indicates the difficulty in optimizing configurations of multiple structural mitigation measures simultaneously, which may lead to a worse objective when running the previous algorithms with more measures.
	Note, however, that the sequential approach generally did not provide improvements over DE-PL for configurations consisting of three, four, and five walls.
	These results indicate that sequential optimization is most beneficial when the number of structural measures becomes larger (e.g., greater than five).

	Finally, a comparison between DE-D-S and the sequential DE-PL solver (DE-PL-S) is made in Table \ref{tab:icold-stats-table-4}.
	On average, DE-PL-S provided a $24\%$ improvement over its direct counterpart.
	The sequential DE-PL solver was also capable of finding a solution which completely mitigated the flood using a smaller structural budget.
	That is, the direct sequential solver found a totally mitigating solution at $m = 9$, but DE-PL-S accomplished this for $m = 8$.
	Interestingly, however, for $m = 10$, DE-D-S found a totally mitigating solution, whereas DE-PL-S only found a \emph{nearly} mitigating solution.
	This again may be a consequence of the relatively small number of experiments performed.
	Overall, except for small $m$ (i.e., $m = 1$), the pathline-based sequential approach appears highly superior to the direct sequential approach.
	This result indicates that DE-PL-S serves as a good general purpose OFMP solver.

\begin{figure}
		\includegraphics[width=1.0\textwidth]{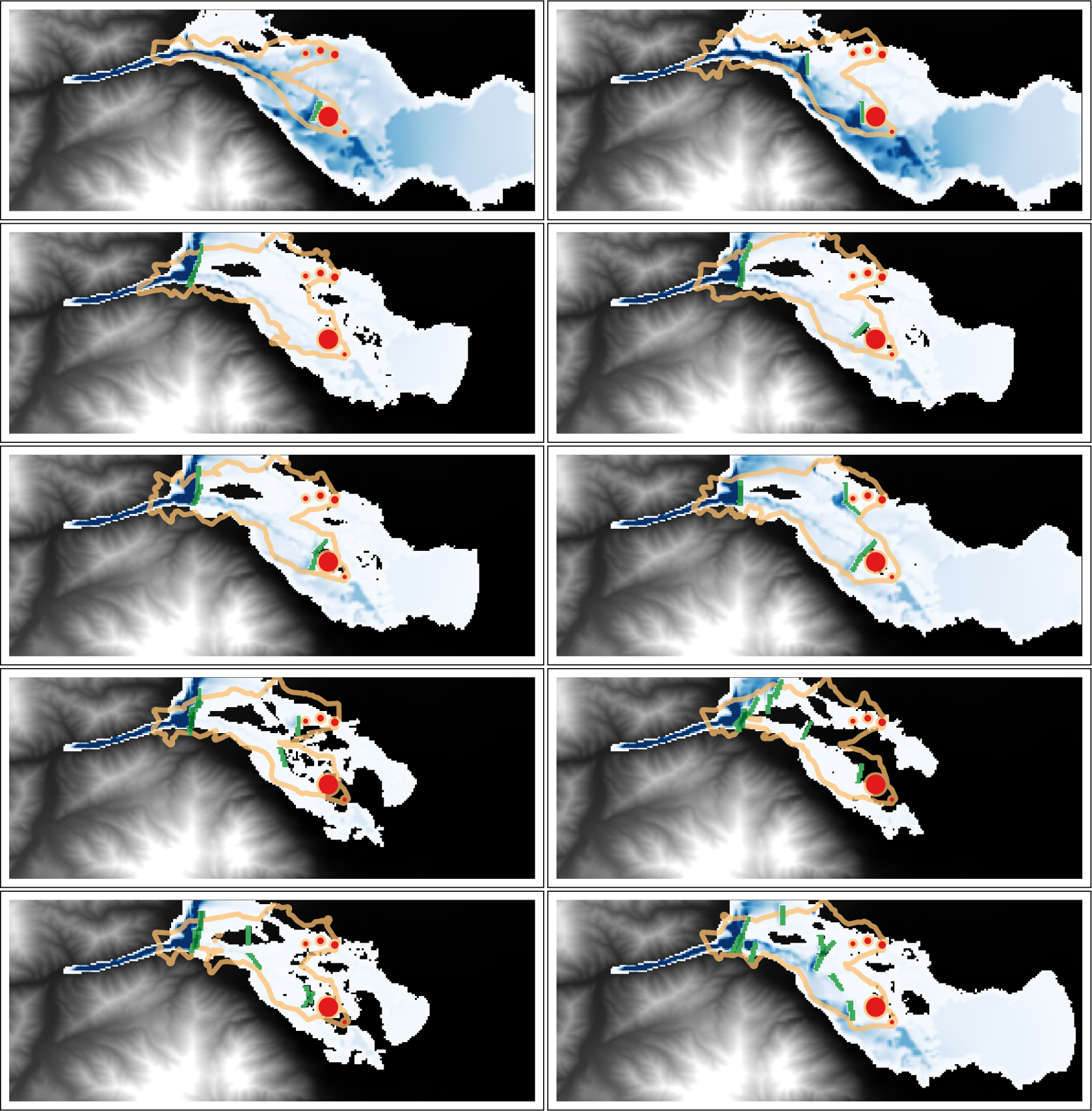}
		\caption{Best obtained elevations and maximum depths for configurations of one through ten walls for the ICOLD 2013 dam failure scenario using DE-PL. Darker blue corresponds to larger maximum depths; gray corresponds to the initial topographic elevation field; green corresponds to elevation additions via the placement of walls; and red corresponds to asset locations. Orange lines represent the exteriors of the final restricted positional sets $\tilde{\mathcal{P}}$ in Algorithm \ref{alg:pseudoalg}.}
		\label{fig:best-icold-summaries}
	\end{figure}

\begin{figure}
		\center
		\includegraphics[width=0.53\textwidth]{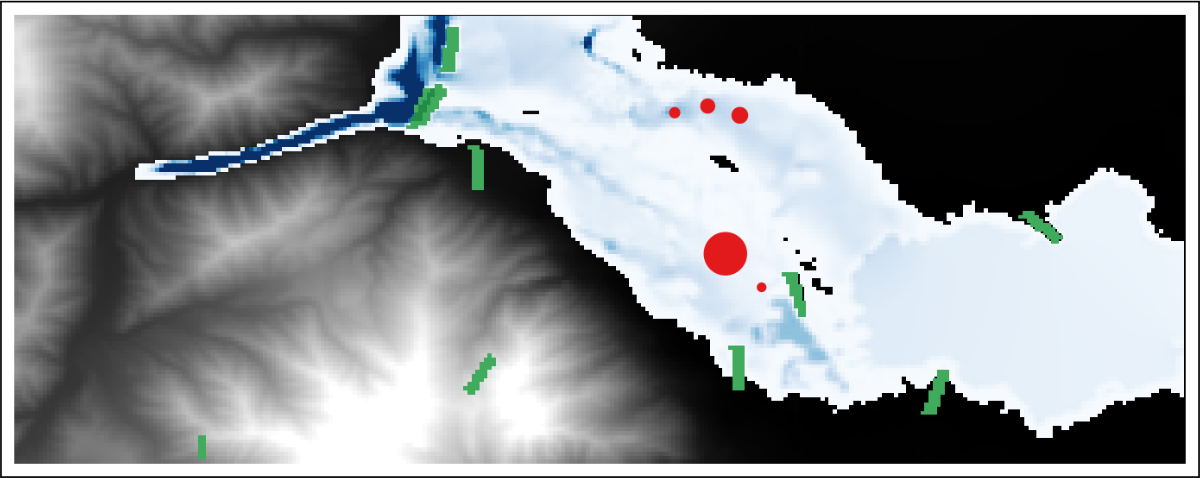}
		\caption{Best solution in a setting equivalent to Figure \ref{fig:best-icold-summaries} for ten walls using DE-D. \hfill}
		\label{fig:de-icold-10}
	\end{figure}

\begin{figure}
		\includegraphics[width=\textwidth]{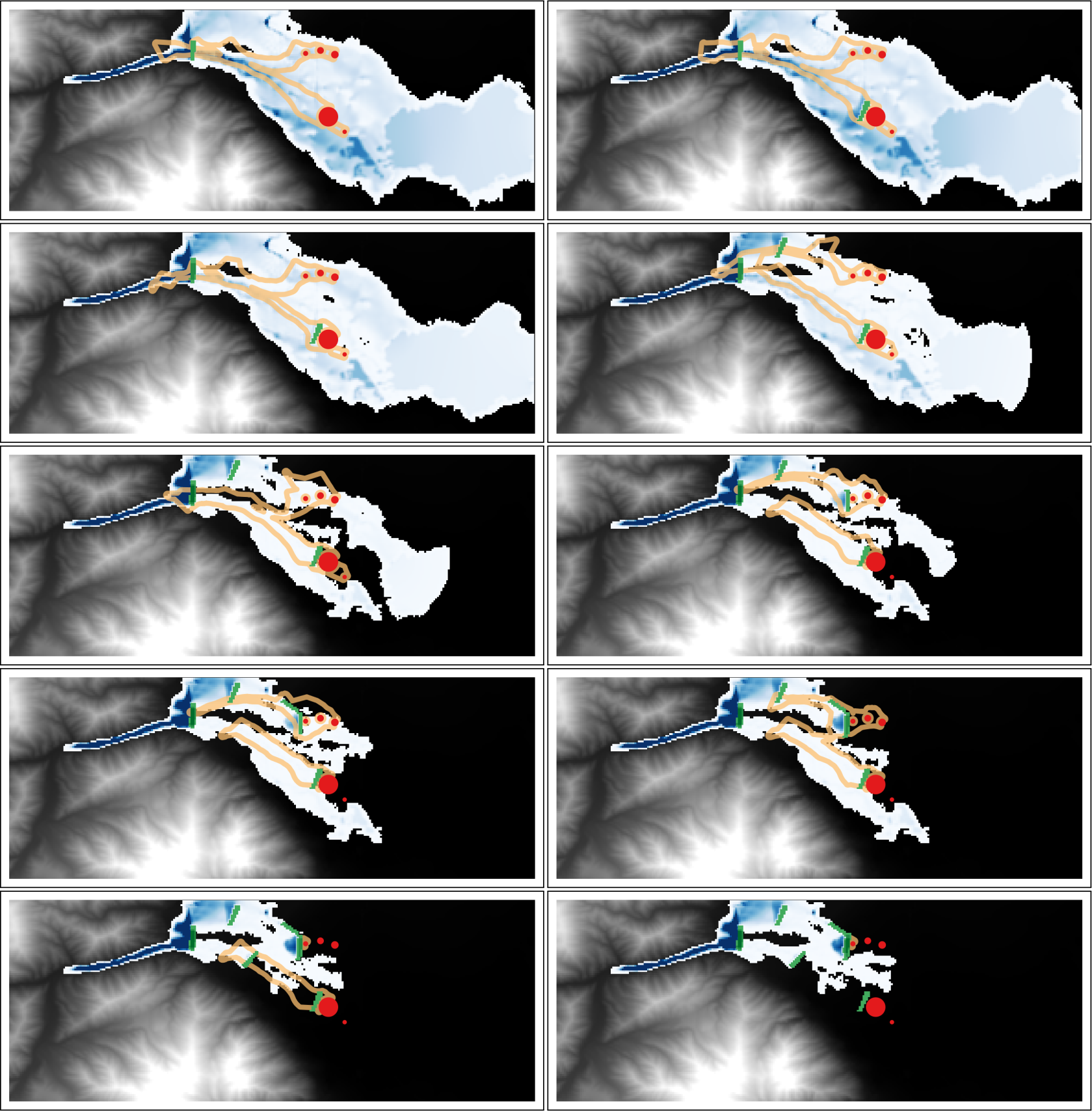}
		\caption{Best obtained elevations and maximum depths for configurations of one through ten walls for the ICOLD 2013 dam failure scenario using DE-PL-S. Darker blue corresponds to larger maximum depths; gray corresponds to the initial topographic elevation field; green corresponds to elevation additions via the placement of walls; and red corresponds to asset locations. Orange lines represent the exteriors of the restricted positional sets $\tilde{\mathcal{P}}$ initialized in Algorithm \ref{alg:pseudoalg}.}
		\label{fig:sequential-icold-summaries}
	\end{figure}

\begin{figure}
		\center
		\includegraphics[width=0.53\textwidth]{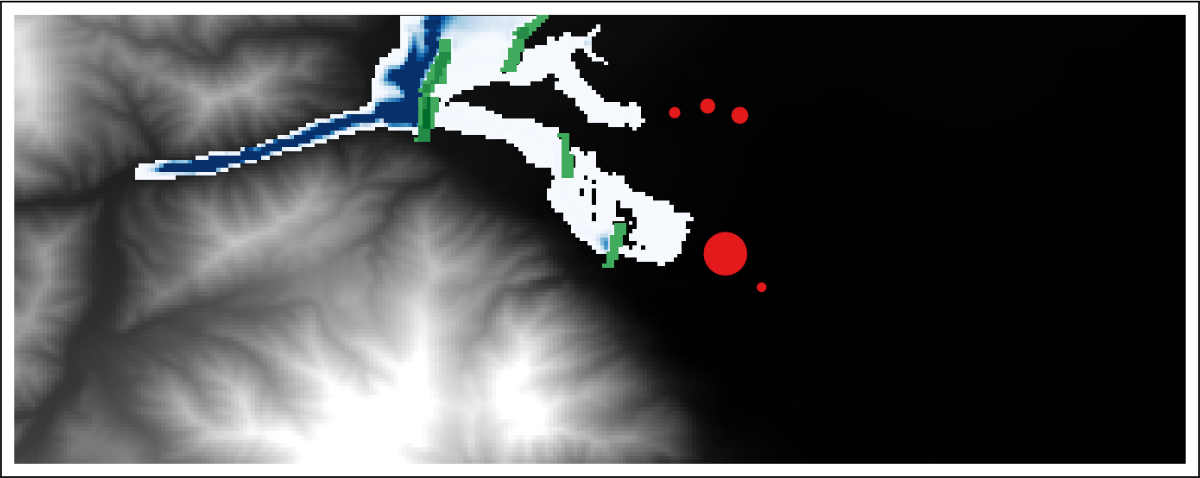}
		\caption{Best solution in a setting equivalent to Figure \ref{fig:sequential-icold-summaries} for ten walls using DE-D-S.}
		\label{fig:sequential-de-icold-10}
	\end{figure}

	Figure \ref{fig:sequential-icold-summaries} displays the ten incremental configurations obtained via DE-PL-S for $m = 10$ and the random seed that gave the minimum corresponding objective in Table \ref{tab:icold-stats-table-4}.
	The ultimate solution for $m = 10$ shows remarkable similarity to the solution obtained via DE-PL for $m = 8$, as shown in Figure \ref{fig:best-icold-summaries}.
	That is, both solutions appear to exploit the critical depth boundary condition to divert water outside of the domain's uppermost boundary.
	However, the sequential solution appears to place a larger number of walls in more intuitive locations.
	Similarly, as displayed by the solution for $m = 10$ shown in Figure \ref{fig:sequential-de-icold-10}, DE-D-S also produced a configuration which diverted flow out of the domain's uppermost boundary, although one wall was placed extraordinarily near this boundary.
	Such solutions may not be possible when using the pathline approach, as pathlines typically do not reside near domain boundaries.

	\subsubsection{Summary of algorithm comparisons}
	\label{sec:alg-summary}
	Tables \ref{tab:icold-stats-table-1} through \ref{tab:icold-stats-table-4} compare the performance of solvers against one another.
	Within these tables, the best objectives over \emph{all} seeds and solvers are denoted in bold, while the best \emph{mean} objectives are underlined.
	It is first apparent that for $m \in \{1, 2\}$, minima were obtained through use of DE-D.
	Good mean objectives were also obtained using this solver.
	This result indicates that direct local search algorithms are capable of performing well on OFMPs that contain a small number of structural measures.
	It also implies that more careful tuning of these algorithms may hold great promise.

	For $m = 3$, DE-PL performed most favorably, providing the best overall and best mean objectives.
	This implies for a moderate number of structural measures, DE-PL effectively uses pathlines to restrict the search space.
	Moreover, if the optimal solution is nonincremental, it is capable of finding solutions that sequential approaches cannot.
	However, for $m > 3$, DE-PL-S performs most favorably, indicating a combination of pathline-based and sequential approaches are needed to solve challenging problems.

	\subsubsection{Proof of concept for soft structural mitigation measures}
	\label{sec:proof-of-concept}
	\begin{figure}[t]
		\includegraphics[width=\textwidth]{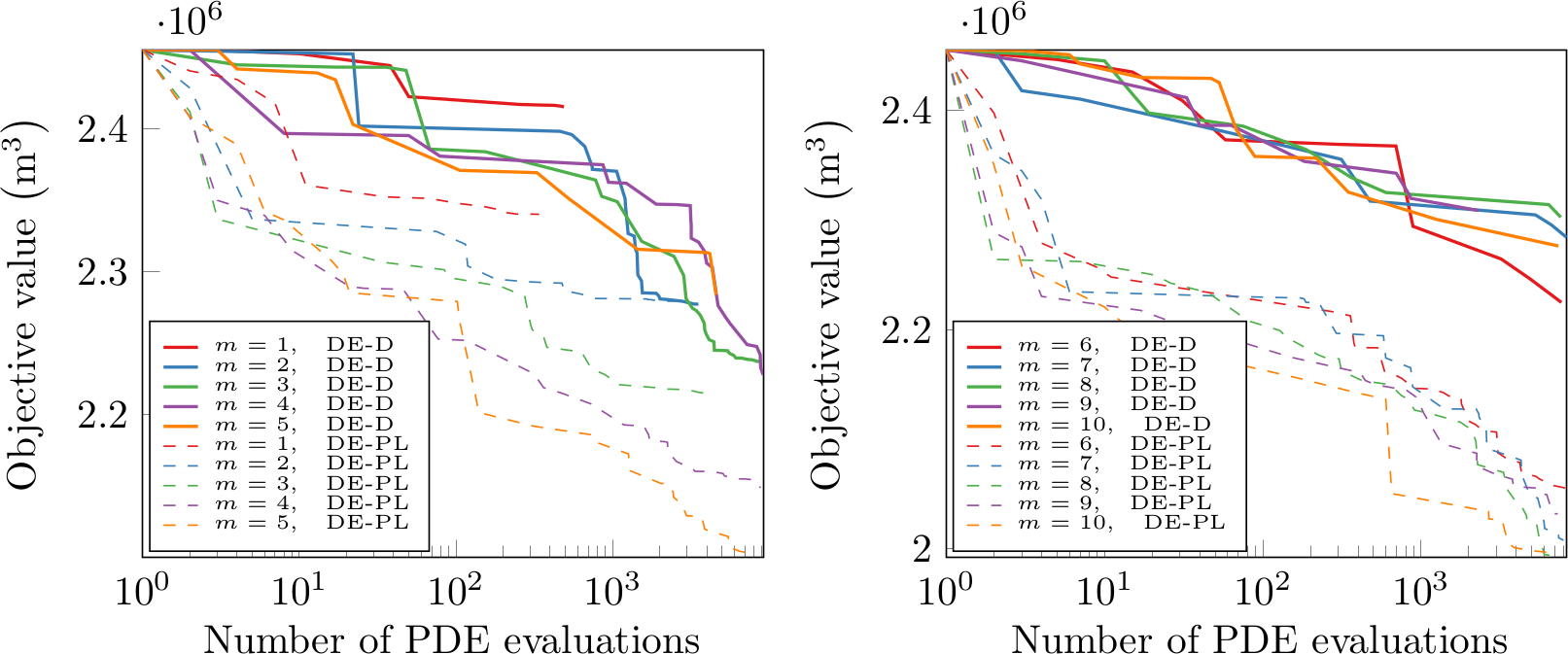}
		\caption{Comparison of objective versus number of PDE evaluations for the OFMP in Section \ref{sec:proof-of-concept}, using DE-D and DE-PL for configurations of one through ten revegetation projects.}
		\label{fig:vegetation-objective-plots}
		\vspace{-3em}
	\end{figure}

	Sections \ref{subsect:results-simple} through \ref{sec:alg-summary} focus on OFMPs designed to configure the placement of \emph{hard} structural mitigation measures (i.e., $\bar{b}_{i} > 0$ and $\bar{n}_{i} = 0$).
	However, it is important to emphasize that the problem formulations and techniques described throughout Sections \ref{sec:model}, \ref{sec:algorithm}, and \ref{subsec:experimental_setting} are not limited to such measures.
	To exemplify this, a proof of concept employing only soft structural measures is assessed.
	In particular, an OFMP taking the form of Equations \eqref{eqn:specific-ofmp-objective-mod} through \eqref{eqn:specific-ofmp-valid-structures-mod} is proposed that optimizes the configuration of $m$ revegetation projects (i.e., $\bar{n}_{i} > 0$ and $\bar{b}_{i} = 0$).

	Using the ICOLD 2013 scenario, the above problem was constructed for a number of revegetation projects ranging from one to ten.
	Each revegetation project was assumed to have a radius of 250 meters and increased the Manning's roughness coefficient in the project region from $0.035$ to $0.123$.
	An experimental setting equivalent to that described in Section \ref{sec:experimental-setting} was used.
	However, in these experiments, only the DE and DE-PL solvers were compared.
	Furthermore, only a single random seed was used.

	In Figure \ref{fig:vegetation-objective-plots}, for each experiment, the objective behavior is plotted against the number of PDE evaluations required to reach that objective.
	The DE-PL solver was generally able to improve upon solutions more rapidly, especially for configurations involving larger numbers of revegetation projects.
	These results mimic the behaviors of Figure \ref{fig:simple-objective-plots}, Table \ref{tab:icold-stats-table-1}, and Table \ref{tab:icold-stats-table-2}.
	That is, for smaller numbers of projects, the direct algorithm is sufficient, but for larger numbers of projects, the pathline-based algorithm is needed to obtain meaningful solutions.

\begin{figure}[!ht]
		\includegraphics[width=\textwidth]{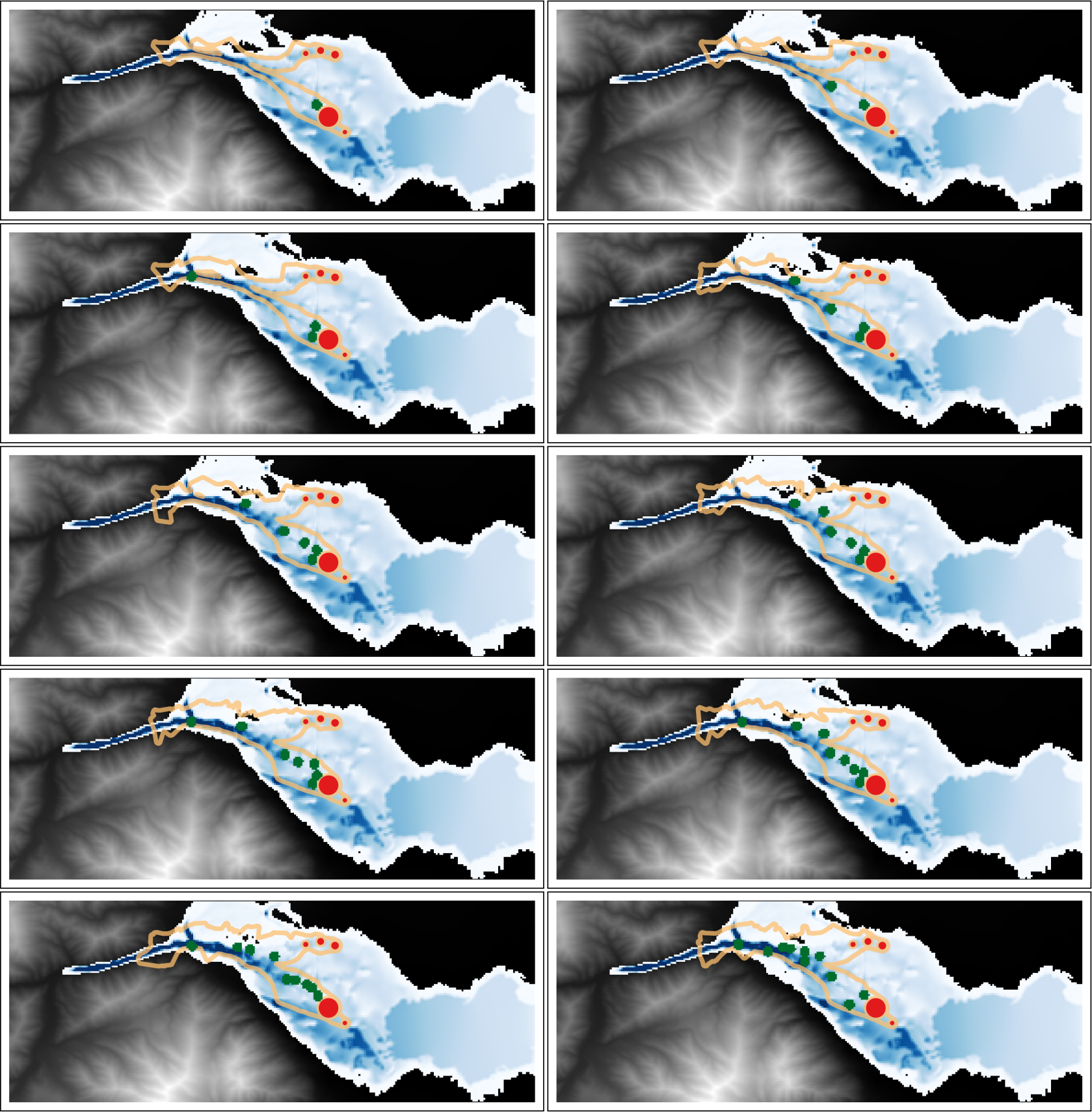}
		\caption{Revegetation locations and maximum depths for configurations of one through ten projects for the ICOLD 2013 dam failure scenario using DE-PL. Darker blue corresponds to larger maximum depths; gray corresponds to the initial topographic elevation field; green corresponds to the placement of revegetation projects; and red corresponds to asset locations. Orange lines represent the exteriors of the restricted positional sets $\tilde{\mathcal{P}}$ initialized in Algorithm \ref{alg:pseudoalg}.}
		\label{fig:best-vegetation-summaries}
	\end{figure}

	Finally, in Figure \ref{fig:best-vegetation-summaries}, the configurations using DE-PL are displayed pictorially for all pairs of scenarios and numbers of projects.
	The results are highly intuitive upon greater inspection.
	First, many of the projects appear to be placed in locations that interdict the initial flood wave.
	More interestingly, many are located along the primary channels of the scenario domain, where larger velocities would occur.
	This makes sense, as the bed shear stress source terms are proportional to the square of velocity; measures that increase roughness are thus highly beneficial in these regions.

	\section{Conclusion}
	\label{sec:conclusion}
	This study addressed the difficult problem of designing structural flood risk management strategies for use within the risk assessment process.
	To this end, an optimization-based decision support approach was proposed for designing mitigation strategies.
	A number of numerical methodologies were developed that generally function through modifying the bed slope and bed shear stress source terms of the 2D shallow water equations.
	However, the methodologies are sufficiently general to modify other source terms (e.g., adjustment of soil properties that affect $\mathbf{S}_{R}$ via infiltration) or even supplant the shallow water equations with a different physical model.

	To formalize the mitigation task, the Optimal Flood Mitigation Problem (OFMP) was introduced.
	To solve practical problems of this type, a time-limited search-based optimization algorithm was developed.
	Within this algorithm, three approaches to generate solutions were explored: a direct approach using only derivative-free optimization, an augmented approach using pathlines to restrict the search space, and a sequential optimization approach.
	The latter two were largely successful, depending on the number of mitigation measures defined in the OFMP.
	Overall, the non-sequential and sequential pathline-based differential evolution approaches provided average improvements of $59\%$ and $65\%$ over their direct counterpart, respectively.
	Results illustrate the first meaningful solutions to large-scale optimization problems of this type.

	Future work should seek to increase and prove the applicability of the approach to realistic flood scenarios.
	First, it should seek to generalize the approach by benchmarking performance on a greater number of real-world flood scenarios.
	Second, it should address the inherent uncertainty in flood scenario parameterizations (e.g., topographic elevation, dam breach parameterization, bed friction).
	To this end, a stochastic optimization approach should be developed to ensure solutions are distributionally robust from a planning perspective.
	Third, a human behavioral study should be conducted to compare the utility of the optimization approach presented herein with the typically manual process used in simulation-based mitigation design.
	Fourth, algorithmic enhancements should be made to increase the realism of mitigation designs. For example, flood walls used in the numerical experiments were overtoppable. This may not be realistic from a flood risk management perspective. Such realism can be embedded within the optimization problem in the form of additional penalties (e.g., when walls are overtopped, a penalty is introduced) or additional constraints.
	Finally, the approach should be extended to solve OFMPs for scenarios that require modeling at finer spatial resolutions.
	To accomplish this, a multi-resolution approach should be developed, where the spatial resolution of a flood scenario is iteratively refined as optimization progresses.
	This work would be valuable for realistic scenarios, where fine resolution details are sometimes necessary to accurately predict flooding behavior.

	\acknowledgments

	All primary input and output data for the experiments described in this study are publicly available at \url{http://nuflood.com/documents/resources/tasseff-2019a.tar.gz}.
	This material is based upon work supported by the National Science Foundation Graduate Research Fellowship under Grant No. DGE-1256260 and under the auspices of the National Nuclear Security Administration of the U.S. Department of Energy at Los Alamos National Laboratory under Contract No. 89233218CNA000001.

	\clearpage


	\clearpage


	\clearpage

	\appendix

	\section{$\textproc{ComputePathline}(\mathbf{U}, x_{0}, y_{0})$}
	\label{sec:pathline-computation}
	\begin{algorithm}[!ht]
		\caption{\textproc{ComputePathline}: Approximates a pathline emanating to some point.}
		 \label{alg:streamline}
		 \begin{algorithmic}[1]
			  \Function{ComputePathline}{$\mathbf{U}, x_{0}, y_{0}$}
					\State $L \gets 0, ~ \ell \gets 0, ~ x \gets x_{0}, ~ y \gets y_{0}, ~ t \gets t_{\text{wet}}(x_{0}, y_{0}), ~ \mathcal{L} \gets \{(x_{0}, y_{0})\}$ \label{alg:streamline:line:initialize}
					\While{$L < L_{\max} ~ \textbf{and} ~ t \in \left[t_{0}, t_{\text{wet}}(x_{0}, y_{0})\right]$} \label{alg:streamline:line:while-loop}
						 \State $(i, j) \gets \textproc{GetIndex}(x, y), ~ k \gets \text{argmin}\left\{\tau \in \mathcal{T}(\mathbf{U}) : |t - \tau|\right\}$ \label{alg:streamline:line:get-current-indices}
						 \If{$\sqrt{u_{ijk}^{2} + v_{ijk}^{2}} \leq \epsilon_{m} ~ \textbf{or} ~ h_{ijk} \leq \epsilon_{m}$} \label{alg:streamline:line:first-check-1}
							  \State \textbf{break} \label{alg:streamline:line:first-check-2}
						 \EndIf \label{alg:streamline:line:first-check-3}
						 \State $\Delta t \gets -\frac{1}{3} \min\left(\frac{\Delta x}{|u_{ijk}|}, ~ \frac{\Delta y}{|v_{ijk}|}\right)$ \label{alg:streamline:line:compute-time-step}
						 \State $x_{*} \gets x + u_{ijk} \Delta t, ~ y_{*} \gets y + v_{ijk} \Delta t, ~ t_{*} = t + \Delta t$ \label{alg:streamline:line:integrate-first}
						 \If{$(x_{*}, y_{*}) \notin D(\mathbf{U})$} \label{alg:streamline:line:second-check-1}
							  \State \textbf{break} \label{alg:streamline:line:second-check-2}
						 \EndIf \label{alg:streamline:line:second-check-3}
						 \State $(i_{*}, j_{*}) \gets \textproc{GetIndex}(x_{*}, y_{*}), ~ k_{*} \gets \text{argmin}\left\{\tau \in \mathcal{T}(\mathbf{U}) : |t_{*} - \tau|\right\}$ \label{alg:streamline:line:get-star-indices}
						 \If{$h_{i_{*}, j_{*}, k_{*}} \leq \epsilon_{m}$} \label{alg:streamline:line:third-check-1}
							  \State \textbf{break} \label{alg:streamline:line:third-check-2}
						 \EndIf \label{alg:streamline:line:third-check-3}
						 \State $x_{n} \gets x + \frac{1}{2} \Delta t \left(u_{ijk} + u_{i_{*}, j_{*}, k_{*}}\right), ~ y_{n} \gets y + \frac{1}{2} \Delta t \left(u_{ijk} + u_{i_{*}, j_{*}, k_{*}}\right)$ \label{alg:streamline:line:integrate-second}
						 \State $\Delta s \gets \sqrt{(x_{n} - x)^{2} + (y_{n} - y)^{2}}$ \label{alg:streamline:line:define-delta-s}
						 \If{$(x_{n}, y_{n}) \notin D(\mathbf{U})$ \textbf{or} $\Delta s \leq \epsilon_{m}$ \textbf{or} $\Delta s > 2 \alpha$} \label{alg:streamline:line:fourth-check-1}
							  \State \textbf{break} \label{alg:streamline:line:fourth-check-2}
						 \EndIf \label{alg:streamline:line:fourth-check-3}
						 \State $L \gets L + \Delta s, ~ \ell \gets \ell + \Delta s$ \label{alg:streamline:line:update-lengths}
						 \State $x \gets x_{n}, ~ y \gets y_{n}, ~ t \gets t_{*}$ \label{alg:streamline:line:final-integration}
						 \If{$\ell \geq \frac{1}{2} (\Delta x + \Delta y)$} \label{alg:streamline:line:save-1}
							  \State $\ell \gets 0, ~ \mathcal{L} \gets \mathcal{L} \cup \{(x_{n}, y_{n})\}$ \label{alg:streamline:line:save-2}
						 \EndIf \label{alg:streamline:line:save-3}
					\EndWhile
					\State \Return $\mathcal{L}$
			  \EndFunction
		 \end{algorithmic}
	\end{algorithm}

	In Algorithm \ref{alg:streamline}, Line \ref{alg:streamline:line:initialize}, the pathline and current pathline segment lengths, $L$ and $\ell$, are initialized to zero, and the pathline-describing point set $\mathcal{L}$ is initialized.
	In Line \ref{alg:streamline:line:while-loop}, the integration loop is defined.
	Integration halts once the total pathline length is greater than some predefined threshold, $L_{\max}$, or the time falls outside the interval of interest, $\left[t_{0}, t_{\text{wet}}\right]$, where $t_{\text{wet}}$ is calculated as per Equation \eqref{eqn:discrete-wet-time}.
	In Line \ref{alg:streamline:line:get-current-indices}, the discrete solution indices are obtained.
	Here, $\textproc{GetIndex}(x, y)$ is a function that maps the spatial coordinates $(x, y)$ to the corresponding spatial index on the rectangular solution grid $G$, $(i, j)$.
	Similarly, the time index $k$ is obtained by computing the index of the ordered timestamp set $\mathcal{T}$ corresponding to the least absolute difference with the current integration time $t$.
	In Lines \ref{alg:streamline:line:first-check-1} through \ref{alg:streamline:line:first-check-3}, the loop is terminated if the current speed or depth is smaller than some arbitrarily small constant $\epsilon_{m}$.

	In Line \ref{alg:streamline:line:compute-time-step}, the time step is computed to (approximately) ensure the integrated distance will not be greater than one third the length of a grid cell.
	In Line \ref{alg:streamline:line:integrate-first}, the first step of second-order Runge-Kutta integration is performed.
	In Lines \ref{alg:streamline:line:second-check-1} through \ref{alg:streamline:line:second-check-3}, the loop is terminated if the point suggested by the previous integration step falls outside the flood scenario's spatial domain, denoted as $D(\mathbf{U})$.
	In Line \ref{alg:streamline:line:get-star-indices}, the discrete indices of the proposed solution are obtained.
	In Lines \ref{alg:streamline:line:third-check-1} through \ref{alg:streamline:line:third-check-3}, the loop is terminated if the depth at the proposed index is too small.
	In Line \ref{alg:streamline:line:integrate-second}, the second Runge-Kutta integration step is performed.
	In Lines \ref{alg:streamline:line:fourth-check-1} through \ref{alg:streamline:line:fourth-check-3}, the loop is terminated if the integrated point falls outside $D(\mathbf{U})$, if the change was small, or if the change was very large (where $\alpha$ is some predefined fixed distance).
	In Line \ref{alg:streamline:line:update-lengths}, the total pathline and temporary segment lengths are updated using the most recent integration distance.
	In Line \ref{alg:streamline:line:final-integration}, the relevant variables are integrated.
	In Lines \ref{alg:streamline:line:save-1} through \ref{alg:streamline:line:save-3}, the temporary segment length is reset to zero and the pathline approximation is updated if the segment length is greater than or equal to the mean grid cell spacing.

	\section{$\textproc{AlphaShape}(Q, \alpha)$}
	\label{sec:alpha-shape-computation}
	In Algorithm \ref{alg:edelsbrunner}, $\textproc{Delaunay}(Q, \alpha)$ is a function that computes the \emph{Delaunay triangulation} for a set $Q$ of discrete points.
	A Delaunay triangulation is a set of triangles such that no point in $Q$ is contained within the circumscribed circle of any triangle.
	A number of algorithms exist to compute this triangulation; herein, that of \citeA{barber1996quickhull} is used.
	\begin{algorithm}[!ht]
		 \caption{\textproc{AlphaShape}: Computes an alpha shape from a discrete set of points $Q$.}
		 \label{alg:edelsbrunner}
		 \begin{algorithmic}[1]
			  \Function{AlphaShape}{$Q, \alpha$}
					\State $\mathcal{D} \gets \textproc{Delaunay}(Q), ~ \mathcal{B} \gets \emptyset$ \label{alg:edelsbrunner:line:initialization}
					\For{$\Delta \in \mathcal{D}$}
						 \State $(a, b, c) \gets \textproc{GetVertices}(\Delta)$ \label{alg:edelsbrunner:line:get-vertices}
						 \State $d_{a} \gets \left\lVert a - b \right\rVert, ~ d_{b} \gets \left\lVert b - c \right\rVert, ~ d_{c} \gets \left\lVert c - a \right\rVert$ \label{alg:edelsbrunner:line:get-edge-distances}
						 \State $s \gets \frac{1}{2}(d_{a} + d_{b} + d_{c})$ \label{alg:edelsbrunner:line:get-semiperimeter}
						 \State $A \gets \sqrt{s \left(s - d_{a}\right) \left(s - d_{b}\right) \left(s - d_{c}\right)}$ \label{alg:edelsbrunner:line:get-area}
						 \If{$A = 0$}
							  \State \textbf{continue} \label{alg:edelsbrunner:line:continue}
						 \ElsIf{$\frac{d_{a} d_{b} d_{c}}{4 A} < \alpha$}
							  \State $\mathcal{B} \gets \mathcal{B} \cup \Delta$ \label{alg:edelsbrunner:line:unionize}
						 \EndIf
					\EndFor
					\State \Return $\mathcal{B}$
			  \EndFunction
		 \end{algorithmic}
	\end{algorithm}

	In Line \ref{alg:edelsbrunner:line:initialization} of Algorithm \ref{alg:edelsbrunner}, the set of Delaunay triangles $\mathcal{D}$ is computed for the point set $Q$, and the set $\mathcal{B}$ comprising the triangular regions of the alpha shape is initialized as the empty set.
	In Line \ref{alg:edelsbrunner:line:get-vertices}, the function $\textproc{GetVertices}(\Delta)$ is used to obtain the vertex positions of the triangle $\Delta$.
	In Line \ref{alg:edelsbrunner:line:get-edge-distances}, the Euclidean edge distances are computed for the triangle $\Delta$.
	In Line \ref{alg:edelsbrunner:line:get-semiperimeter}, the semiperimeter $s$ of the triangle $\Delta$ is computed.
	In Line \ref{alg:edelsbrunner:line:get-area}, the area of the triangle $\Delta$ is computed via Heron's formula.
	In Line \ref{alg:edelsbrunner:line:unionize}, if the circumscribed radius of the triangle is less than the constant $\alpha$, the triangle is unioned with the set $\mathcal{B}$ describing the alpha shape.
	In this paper, $\alpha$ is always taken to be $5 (\Delta x + \Delta y) / 2$, where $\Delta x$ and $\Delta y$ are the grid cell spacings used to discretize the spatial domain in the $x$- and $y$- directions, respectively.
\end{document}